\documentclass[3p, 10pt, number]{elsarticle}
\biboptions{sort}
\newdefinition{remark}{Remark}
\usepackage{amssymb,amsmath,amsfonts}
\usepackage{pgf,pgfplots,pgfplotstable}

\usepackage{upgreek} 

\pgfplotsset{compat=newest}
\usepgfplotslibrary{groupplots, polar, external}
\usetikzlibrary{shapes,arrows,decorations,calendar,matrix,backgrounds,folding,calc,positioning,spy, arrows.meta, bending, patterns}

\tikzsetexternalprefix{./figures/external/}
\tikzset{external/only named=true}
\tikzset{external/mode=only graphics}
\definecolor{matblue}{rgb}{0 0.4470 0.7410}
\definecolor{matorange}{rgb}{0.8500 0.3250 0.0980}
\definecolor{matyellow}{rgb}{0.9290 0.6940 0.1250}
\definecolor{matpurple}{rgb}{0.4940 0.1840 0.5560}
\definecolor{matgreen}{rgb}{0.4660 0.6740 0.1880}
\definecolor{matazure}{rgb}{0.3010 0.7450 0.9330}
\definecolor{matred}{rgb}{0.6350 0.0780 0.1840}
\definecolor{amber}{rgb}{1.0, 0.75, 0.0}

\usepackage{tabularx, booktabs, hyperref}
\newcolumntype{C}{>{\centering\arraybackslash}X}
\newcolumntype{L}{>{\raggedright\arraybackslash}X}
\newcolumntype{s}{>{\hsize=.5\hsize}C}
\newcolumntype{T}{>{\hsize=.05\hsize\raggedright\arraybackslash}X}
\newcolumntype{M}{>{\hsize=.25\hsize}C}

\usepackage[capitalize]{cleveref}  

\usepackage[ruled]{algorithm2e} 
\SetKwInput{KwReq}{Require}
\crefname{algorithm}{algorithm}{algorithms}
\crefname{algocf}{algorithm}{algorithms}
\Crefname{algorithm}{Algorithm}{Algorithms}

\usepackage{lmodern, float, bm, caption, url}

 \usepackage{setspace}

\renewcommand{\hat}{\widehat}
\renewcommand{\tilde}{\widetilde}

\newcommand{\pdt}[1]{{\frac{\partial #1}{\partial t}}}

\newcommand{\diag}{{\operatorname{diag}}}
\newcommand{\tr}{{\operatorname{tr}}}

\newcommand{\grad}{{\nabla}}

\renewcommand{\Re}{{\operatorname{Re}}}

\newcommand{\Ma}{{\operatorname{Ma}}}

\newcommand{\St}{{\operatorname{St}}}

\newcommand{\avg}[1]{\{\!\!\{#1\}\!\!\}}
\newcommand{\jump}[1]{{{[\![#1]\!]}}}

\newcommand{\facets}{\mathcal{F}_h}
\newcommand{\mesh}{\mathcal{T}_h}

\newcommand{\mat}[1]{{\bm{#1}}}
\renewcommand{\vec}[1]{{\bm{#1}}}

\newcommand{\I}{{\mat{I}}}
\newcommand{\T}{{\mathrm{T}}}

\newcommand{\VEL}{{\vec{u}}}

\newcommand{\EPS}{{\mat{\varepsilon}}}

\journal{Computational Physics}

\begin{document}
\graphicspath{{./figures/}}
\begin{frontmatter}

    \title{IMEX Schemes for Compressible Flow using Hybridizable Discontinuous Galerkin Methods} 

    \author[1]{Jan Ellmenreich\corref{cor1}}
    \ead{jan.ellmenreich@tuwien.ac.at}
    \cortext[cor1]{Corresponding author}

    \author[2]{Edmond K. Shehadi}

    \author[2]{Philip L. Lederer}

    \affiliation[1]{organization={Institute of Analysis and Scientific Computing, TU~Wien},
        addressline={Wiedner~Hauptstrasse~8-10},
        city={Vienna},
        postcode={1040},
        country={Austria}}
    \affiliation[2]{organization={Fachbereich Mathematik, Universität~Hamburg},
        addressline={Bundesstraße~55},
        city={Hamburg},
        postcode={20146},
        country={Germany}}

    \begin{abstract}
        In this work, we develop a geometry-split implicit-explicit (IMEX) framework for the compressible flow equations, wherein stiff regions are treated via an implicit hybridizable discontinuous Galerkin (HDG) method, while non-stiff regions are treated via an explicit discontinuous Galerkin (DG) method. Two implicit formulations are investigated: a mixed HDG method (HDG-MX) and a primal interior-penalty HDG method (HDG-IP). The spatial coupling between the implicit and explicit solutions is achieved in a conservative manner by appropriate interface conditions, while the temporal synchronization is maintained through the use of additive Runge-Kutta (ARK) schemes. We provide a detailed discussion on the computational performance of the resulting IMEX schemes. Verification and validation over a range of numerical experiments confirm that the proposed IMEX schemes achieve high-order accuracy in both space and time. Performance studies further indicate that the approach effectively alleviates geometry-induced stiffness and can provide speedups of up to approximately $50$ relative to a fully explicit DG scheme, provided that the implicit region is chosen appropriately.
    \end{abstract}




\end{frontmatter}

\section{Introduction}
\label{sec::introduction}

Stiff problems are ubiquitous in computational fluid dynamics (CFD) and remain a central challenge in the efficient simulation of complex flows. A primary source of stiffness is the disparity of spatial and temporal scales, particularly pronounced in turbulent boundary layers \cite{pope2001turbulent}. This scale separation is especially relevant in wall-resolved large-eddy simulation (LES) and direct numerical simulation (DNS), where accurate near-wall resolution is essential. Consequently, stiffness is often highly localized and associated with anisotropic and distorted mesh elements, leading to so-called \textit{geometry-induced stiffness}.

High-fidelity simulations increasingly rely on high-order spatial discretizations, such as discontinuous Galerkin (DG) \cite{reed1973triangular} and hybridizable discontinuous Galerkin (HDG) \cite{cockburn2009unified} methods, due to their superior accuracy and favorable dispersion and dissipation properties \cite{ainsworth2004106,wang2013high}. However, these discretizations can exacerbate stiffness: the CFL constraint is governed by an \textit{effective} element length scale that decreases with increasing polynomial degree \cite{hesthaven2008nodal,chalmers2020109095}, thereby restricting explicit time integration.

Fully implicit schemes provide a natural mechanism to alleviate stiffness by allowing larger time steps. However, they require the solution of large, globally coupled nonlinear systems at each step, which becomes prohibitively expensive for simulations involving millions or billions of degrees of freedom (DOFs) \cite{slotnick2014cfd}. This limitation is particularly severe in industrially relevant flows at high Reynolds numbers ($\mathrm{Re} > 10^6$), where the computational cost of unsteady simulations scales steeply. For instance, the number of DOFs scales as $N_{\mathrm{DNS}} \sim \mathrm{Re}^{37/14}$ for DNS and $N_{\mathrm{WR\text{-}LES}} \sim \mathrm{Re}^{13/7}$ for wall-resolved LES \cite{choi2012grid}, underscoring the prohibitive cost of fully implicit approaches in such regimes.

Fully explicit schemes, in contrast, avoid global coupling and are attractive due to their simplicity, low memory requirements and excellent parallel scalability. Their applicability, however, is severely limited by stringent stability constraints, particularly in subsonic regimes. In stiff regions--especially near walls--the time step is dictated by the smallest grid spacing, rendering explicit methods inefficient for complex geometries and long-time integrations. Although local time-stepping (LTS) techniques \cite{dumbser20083971,zanotti2015204,choi2012grid,breuer2016petascale} can partially alleviate these restrictions, their effectiveness in large-scale simulations is often limited by load imbalance and increased algorithmic complexity.

Implicit-explicit (IMEX) time integration strategies offer a natural compromise. These methods can be broadly classified into \textit{semi-implicit} and \textit{domain-partitioned} approaches, which target different sources of stiffness. Semi-implicit methods treat fast processes (e.g., diffusive or acoustic effects) implicitly to address scale-separated stiffness \cite{ascher1997implicit}, whereas domain-partitioned approaches target geometry-induced stiffness by applying implicit integration only in localized regions \cite{kanevsky20071753}.

Following the seminal work of \cite{ascher1997implicit}, Runge-Kutta (RK)-based IMEX schemes have become a widely adopted class of methods for stiff systems. Their popularity stems from their ability to achieve high-order temporal accuracy within a multi-stage, one-step framework. Typically, an explicit RK method is coupled with a singly diagonally implicit RK (SDIRK) scheme for the stiff components \cite{kennedy2016diagonally,kennedy2019221}. The SDIRK structure enables stage-by-stage solves with a constant diagonal coefficient, simplifying the implicit treatment and making it well suited for large-scale applications, particularly when combined with quasi-Newton methods.

These IMEX formulations have been successfully combined with high-order discontinuous discretizations. Within DG frameworks, \cite{kanevsky20071753} addressed geometry-induced stiffness in compressible flow, while \cite{kang2020109010} proposed a hybrid HDG-DG approach for shallow water equations in the presence of scale-separated stiffness. More recently, \cite{vermeire2021110022} introduced an accelerated IMEX (AIMEX) strategy within a flux reconstruction (FR) framework, achieving noticeable speedups for compressible flows affected by geometry-induced stiffness. Their approach introduces additional intermediate stages that enhance the stability of the explicit terms without increasing the number of implicit solves. This approach was later extended in \cite{pereira2025113819} using hybridized and embedded FR variants, leading to more substantial efficiency gains.

In this work, we adopt a hybrid IMEX strategy that combines DG discretizations for the explicit terms with HDG discretizations for the implicit terms, yielding a more efficient alternative to the classical DG-DG IMEX formulation of \cite{kanevsky20071753}. While related hybrid approaches have been proposed \cite{kang2020109010,pereira2025113819}, existing studies either focus on shallow water equations or FR-based discretizations for compressible flows, with limited attention devoted to the treatment of viscous terms. This is a key limitation, as viscous terms are well known to dominate the computational cost in (H)DG methods due to gradient evaluations.

Motivated by this observation, we investigate efficient viscous discretizations within an HDG framework that fully exploit the RK-IMEX setting and the proposed HDG-DG splitting for geometry-induced stiffness. In particular, we employ a standard DG interior penalty (DG-IP) discretization for the explicit terms and consider two alternative HDG formulations for the implicit operator. The first is a mixed formulation (HDG-MX), which introduces auxiliary variables for the symmetric rate-of-strain tensor and temperature gradient \cite{vila-perezHybridisableDiscontinuousGalerkin2021,ellmenreichCharacteristicBoundaryConditions2026}. The second is an interior penalty formulation (HDG-IP) that avoids such auxiliary variables entirely.

While \textit{primal} HDG formulations for compressible flows have been explored previously \cite{dahm2017toward,franciolini2020104542}, they rely on BR2-type viscous discretizations \cite{bassi1997high}, which incur non-trivial computational overhead due to lifting operations. In contrast, the HDG-IP formulation considered here leads to a more compact and computationally efficient implementation, particularly for element-local operations. To the best of our knowledge, HDG formulations based on interior penalty discretizations remain largely unexplored in compressible flow simulations, despite their favorable computational properties--an aspect that becomes increasingly important as the cost of local operations grows with the number of DOFs, especially in three dimensions \cite{dahm2017toward}.

The contributions of this work are threefold. First, we develop an HDG-DG IMEX framework for the compressible Navier-Stokes equations. Second, we provide a detailed comparison between HDG-MX and HDG-IP formulations, highlighting the trade-offs between formulation complexity and computational efficiency. Third, we demonstrate that combining IMEX time integration with an explicit DG-IP discretization and an implicit HDG-IP discretization provides an effective strategy for mitigating geometry-induced stiffness in compressible flow simulations. The proposed methodology is implemented within the open-source \href{https://plederer.github.io/dream_solver/index.html}{\texttt{DreAm}} package\footnote{\url{https://plederer.github.io/dream_solver/index.html}}.

The remainder of this paper is organized as follows. In \Cref{sec::compressible_equations}, we introduce the governing equations of fluid flow. The spatial discretization is then presented in \Cref{sec::spatial_discretization}, where the DG-IP, HDG-IP and HDG-MX formulations are described in \Cref{sec::sdg,sec::hdg_ip,sec::hdg_mixed}, respectively. The IMEX time integration strategy, along with the static condensation procedure for the implicit system and the theoretical computational performance, is discussed subsequently in \Cref{sec::imex}. Numerical results are presented in \Cref{sec::experiments}, followed by concluding remarks in \Cref{sec::conclusion}.

\section{Compressible flow equations}
\label{sec::compressible_equations}

The two-dimensional compressible Navier-Stokes equations in conservative form are
\begin{align}
	\label{eq::compressible_equations}
    \frac{\partial \vec{U}}{\partial t}
    +
    \nabla \cdot
    \left( 
    \vec{F}^{c} - \vec{F}^{\nu}
    \right)
    =
    \vec{0},
\end{align}
where the vector of conservative variables is $\vec{U} = (\rho, \rho u, \rho v, \rho E)^T$. Here, $\rho$ is the density, $(u, v)$ are the velocity components and $E$ is the total energy per unit mass, defined as $E = e_i + e_k$, where the internal energy is $e_i = c_v \theta$ (assuming a calorically perfect gas) and kinetic energy is $e_k = (u^2 + v^2) / 2$, both expressed per unit mass. Here and in the following, repeated indices imply summation (Einstein notation).

The inviscid (convective) and viscous fluxes are given by
\begin{align} \label{eq::inviscid_viscous_fluxes}
	\vec{F}^{c}
	=
	\begin{pmatrix}
		\vec{F}^{c}_1,\,
		\vec{F}^{c}_2
	\end{pmatrix}
	=
	\begin{pmatrix}
		\rho u        & \rho v\\
		\rho u^2 + p  & \rho u v\\
		\rho u v      & \rho v^2 + p\\
		\rho u H      & \rho v H
	\end{pmatrix},
	\qquad
	\vec{F}^{\nu}
	=
	\begin{pmatrix}
		\vec{F}^{\nu}_1,\,
		\vec{F}^{\nu}_2
	\end{pmatrix}
	=
	\begin{pmatrix}
		0                       &  0\\
		\tau_{11}               &  \tau_{12}\\
		\tau_{21}               &  \tau_{22}\\
		\tau_{1j} u_{j} - q_{1} &  \tau_{2j} u_{j} - q_{2}
	\end{pmatrix},
\end{align}
where $\vec{F}_i^c$ and $\vec{F}_i^\nu$ denote the fluxes in the $x_i$-direction ($i=1,2$), $H = E + p/\rho$ is the total enthalpy per unit mass, $\tau_{ij}$ is the viscous stress tensor and $q_i$ is the heat flux vector.

For subsequent developments, it is convenient to express the fluxes in terms in their quasi-linear form
\begin{align} \label{eq::definition_jacobians_inviscid_viscous}
	\vec{F}^{c}_{i} (\vec{U})
	=
	\mat{A}_i (\vec{U})\, \vec{U},
	\qquad
	\vec{F}_i^\nu (\vec{U},\, \nabla \vec{U}) = \mat{K}_{ij}(\vec{U})\, \frac{\partial \vec{U}}{\partial x_j},
\end{align}
where $\mat{A}_{i} = \partial \vec{F}^{c}_{i} / \partial \vec{U}$ and $\mat{K}_{ij} = \partial \vec{F}^{\nu}_{i} / \partial (\partial \vec{U} / \partial x_j)$ are the inviscid and viscous Jacobian matrices, respectively, defined in \ref{app::jacobians}. Accordingly, we introduce the compact notation
\begin{align} \label{eq::diffusion_block_matrix_product}
	\mat{K} \nabla \vec{U} := \left( \vec{F}_1^\nu,\, \vec{F}_2^\nu \right).
\end{align}

The pressure is related to the thermodynamic variables through the ideal gas equation of state
\begin{align}
    p = \rho R \theta,
\end{align}
where $\theta$ is the temperature and $R = c_p - c_v$ the specific gas constant, with $c_p$ and $c_v$ denoting the specific heats at constant pressure and constant volume, respectively. The heat flux is modeled using Fourier's law
\begin{align}
    q_{i} = - \kappa \frac{\partial \theta}{\partial x_{i}},
\end{align}
where $\kappa$ is the thermal conductivity, taken constant: $\kappa = c_p \mu / \mathrm{Pr}$. In this work, $\mathrm{Pr} = 0.72$ is chosen and the dynamic viscosity $\mu$ is held constant.

The viscous stress tensor for a Newtonian fluid is given by
\begin{align}
    \tau_{ij}
    =
    \mu \left( \frac{\partial u_j}{\partial x_i} + \frac{\partial u_i}{\partial x_j} \right)
    +
    \lambda \frac{\partial u_k}{\partial x_k} \delta_{ij},
\end{align}
where $\delta_{ij}$ is the Kronecker delta and $\lambda = -2\mu/3$ is the second viscosity from Stokes' hypothesis.

\section{Spatial Discretization}
\label{sec::spatial_discretization}

Three variants of the modal discontinuous Galerkin (DG) method are considered: a standard DG formulation and two hybridizable DG (HDG) formulations. Viscous terms are treated using an interior penalty (IP) approach in DG, while HDG is implemented in both IP and mixed forms. In all cases, convective terms are discretized using an \textit{upwind} flux-vector splitting strategy.

\subsection{Mesh partition and functional spaces}
Consider a mesh $\mesh$ to be the partition of a domain $\Omega$ into shape-regular, uniform and non-overlapping polyhedral mesh elements $T$. The mesh is assumed to satisfy the following regularity conditions:
\begin{align*}
	\max_{T \in \mesh} \frac{h_T^d}{|T|} & \leq C    \qquad C > 0 \quad \forall T \in \mesh, & \qquad
	h_T                                  & \leq C h  \qquad C > 0 \quad \forall T \in \mesh,
\end{align*}
where $h_T$ denotes the element diameter, and the global mesh size $h$ is defined as
\begin{align*}
	h_T  := \text{diam}(T), \qquad h  := \max\limits_{T \in \mesh}{h_T}.
\end{align*}

Let $\facets$ denote the set containing the mesh skeleton,
which consists of the sets of interior facets $\facets^{\text{int}}$ and exterior facets $\facets^{\text{ext}}$, given by
\begin{align*}
	\facets              := \bigcup_{T \in \mesh} \partial T,   \qquad
	\facets^{\text{int}} := \facets \backslash \partial \Omega, \qquad
	\facets^{\text{ext}} := \facets \cap \partial \Omega.
\end{align*}

On an element $T$, the outward unit normal vector \cite{dipietroMathematicalAspectsDiscontinuous2012} on $\partial T$ is denoted by $\vec{n}$. On a common facet $F = \partial T^+ \cap \partial T^- \in \facets^{\text{int}}$, the outward unit normal vectors with respect to the neighboring elements $T^+$~and~$T^-$ are denoted by $\vec{n}^+$ and $\vec{n}^-$, respectively, such that $\vec{n}^+ = -\vec{n}^-$.

Denote by $L^2(\Omega)$ the space of real-valued square-integrable functions over $\Omega$, equipped with the inner product and norm
\begin{align*}
	\left( u, v \right) := \sum_{T \in \mesh} \int_{T} u v \, d\bm{x}, \qquad
	\|u\|               := \sqrt{( u, u )},
\end{align*}
and denote by $L^2(\facets)$ the space of real-valued square-integrable functions on $\facets$, with inner product and norm
\begin{align*}
	{\left\langle u, v \right\rangle}_{\facets} := \sum_{F \in \facets} \int_{F} u v \, d\bm{s}, \qquad
	\|u\|_{\facets}                             := \sqrt{\left\langle u, u \right\rangle}_{\facets}.
\end{align*}
For an arbitrary space $V$, define $V(\Omega, \mathbb{R}^n)$ as the vector-valued space, whose $n$-components are in~$V$, and $V(\Omega, \mathbb{R}^{n \times n})$
as the matrix-valued space, whose $n^2$-components are in $V$. Then $L^2(\Omega, \mathbb{R}^{n \times n}_{\mathrm{sym}})$ denotes the space of 
real-valued square-integrable symmetric tensors over $\Omega$, defined as
\begin{align*}
	L^2(\Omega, \mathbb{R}^{n \times n}_{\mathrm{sym}}) := \left\{ \mat{v} \in L^2(\Omega, \mathbb{R}^{n \times n}) : \mat{v} = \mat{v}^\T \right\}.
\end{align*}

The spaces of polynomials of degree at most $k$, restricted to an element $T$ or a facet $F$, are denoted by $\mathbb{P}^k(T)$ and $\mathbb{P}^k(F)$, respectively.
The broken polynomial spaces on $\mesh$ and $\facets$ are then given by:
\begin{align*}
	\mathbb{P}^k(\mesh)                 & := \lbrace v \in L^2(\Omega) : v|_T \in \mathbb{P}^k(T) \, \forall T \in \mesh \rbrace,                                            \\
	\mathbb{P}^k(\facets)               & := \lbrace v \in L^2(\facets) : v|_F \in \mathbb{P}^k(F) \, \forall F \in \facets \rbrace,                                         \\
	\mathbb{P}^k(\mesh, \mathbb{R}^n)   & := \lbrace \vec{v} \in L^2(\Omega, \mathbb{R}^n) : \vec{v}|_T \in \mathbb{P}^k(T, \mathbb{R}^n) \, \forall T \in \mesh \rbrace,    \\
	\mathbb{P}^k(\facets, \mathbb{R}^n) & := \lbrace \vec{v} \in L^2(\facets, \mathbb{R}^n) : \vec{v}|_F \in \mathbb{P}^k(F, \mathbb{R}^n) \, \forall F \in \facets \rbrace, \\
	\mathbb{P}^k(\mesh, \mathbb{R}^{n \times n}_{\mathrm{sym}})   & := \lbrace \vec{v} \in L^2(\Omega, \mathbb{R}^{n \times n}_{\mathrm{sym}}) : \vec{v}|_T \in \mathbb{P}^k(T, \mathbb{R}^{n \times n}_{\mathrm{sym}}) \, \forall T \in \mesh \rbrace.
\end{align*}

On a given time interval $(0, t_{\textrm{total}}]$, we define the following finite element spaces
\begin{subequations}
	\label{eq::finite_element_spaces}
	\begin{alignat}{1}
    U_h       & := L^2\left( (0, t_{\textrm{total}}] ; \mathbb{P}^k(\mesh, \mathbb{R}^{4})          \right),      \\
    \hat{U}_h & := L^2\left( (0, t_{\textrm{total}}] ; \mathbb{P}^k(\facets, \mathbb{R}^{4})         \right),              \\
	Q_h 	 & := L^2\left( (0, t_{\textrm{total}}] ; \mathbb{P}^k(\mesh, \mathbb{R}^{2 \times 2}_{\mathrm{sym}}) \right) \times  L^2\left( (0, t_{\textrm{total}}] ; \mathbb{P}^k(\mesh, \mathbb{R}^{2}) \right).
\end{alignat}
\end{subequations}

\subsection{Characteristic Flux Splitting}
\label{subsec:inviscid_flux_splitting}

Using \Cref{eq::definition_jacobians_inviscid_viscous}, the Jacobian of the inviscid flux projected along the outward unit normal $\vec{n}$ admits the eigendecomposition
\begin{align}
	\mat{A}_n (\vec{U})
	:=
	\mat{A}_i (\vec{U}) n_i
	=
	\mat{P} (\vec{U})\, \mat{\Lambda} (\vec{U})\, \mat{P}^{-1} (\vec{U}),
\end{align}
where $\mat{P}$ is the matrix of right eigenvectors and $\mat{\Lambda}$ is the diagonal matrix of eigenvalues associated with the normal direction.

The Jacobian can be decomposed into outgoing and incoming contributions according to the sign of the eigenvalues,
\begin{align}\label{eq::inviscid_jacobian_eigendecomposition}
	\mat{A}_n^{\pm} (\vec{U})
	=
	\mat{P} (\vec{U})\,\mat{\Lambda}^{\pm} (\vec{U})\, \mat{P}^{-1} (\vec{U}),
	\qquad
	\mat{\Lambda} (\vec{U}) = \mat{\Lambda}^{+} (\vec{U}) + \mat{\Lambda}^{-} (\vec{U}),
\end{align}
where $\mat{\Lambda}^{+}$ and $\mat{\Lambda}^{-}$ are diagonal matrices containing the positive (outgoing) and negative (incoming) eigenvalues of $\mat{\Lambda}$, respectively, with zeros elsewhere.

\subsection{Standard Discontinuous Galerkin: Interior Penalty}
\label{sec::sdg}

For the DG, the inviscid fluxes are treated using a standard Riemann solver at element surfaces, while an adjoint-consistent symmetric interior penalty (IP) formulation \cite{hartmann2008optimal,hartmann2007adjoint} is used to discretize the viscous terms. 

The resulting weak formulation of the DG-IP scheme reads: find $\vec{U}_{h} \in U_h$, such that
\begin{align} \label{eq::weak_form_sdg_ip}
	\sum_{T \in \mesh} 
	\int\limits_{T} \pdt{\vec{U}_{h}} \cdot \vec{V}_{h} \, d\bm{x}
	+
	\mathcal{B}_{I} \big( \vec{U}_{h}; \vec{V}_{h} \big) 
	+
	\mathcal{B}_{\Gamma} = 0,
	\qquad 
	\forall\ \vec{V}_{h} \in U_h,
\end{align}
where $\mathcal{B}_{I}$ and $\mathcal{B}_{\Gamma}$ are forms emanating from the spatial discretization of the fluxes in the interior and exterior of the domain, respectively. 

In this section, the average and jump operators for the DG are defined as
\begin{align}
	\avg{\vec{\psi}} = \frac{\vec{\psi}^{+} + \vec{\psi}^{-}}{2}, 
	\qquad
	\jump{\vec{\psi}} = \vec{\psi}^{+} \otimes \vec{n}^{+} + \vec{\psi}^{-} \otimes \vec{n}^{-},
\end{align}
where superscripts $(\cdot)^{-}$ and $(\cdot)^{+}$ correspond to the two neighboring solutions sharing a common face, with the notation $\vec{\psi} \otimes \vec{n} = \psi_i n_j$, implying the usual outer product written in Einstein/tensor notation.

In \Cref{eq::weak_form_sdg_ip}, boundary conditions are imposed through the operator $\mathcal{B}_{\Gamma}$, whose definition is provided in \ref{sec::sdg_bcs}. In the present work, we restrict attention to far-field and adiabatic wall boundary conditions. Additional formulations, particularly those addressing non-reflecting boundary conditions within the DG-IP framework, can be found in \cite{shehadi2024109194,shehadi2024NRBC}. The remaining form, $\mathcal{B}_{I}$, is defined as
\begin{multline} \label{eq::bilinear_form_interior_sdg}
	\mathcal{B}_{I} \big( \vec{U}_{h}; \vec{V}_{h} \big)
	=
	\sum_{T \in \mesh} 
	\int\limits_{T}
	\Big(
	\vec{F}^{\nu} (\vec{U}_{h}, \nabla \vec{U}_{h})
	-
	\vec{F}^{c}(\vec{U}_{h}) 
	\Big) : \grad{\vec{V}_{h}} \, d\bm{x}
	+
	\hspace{-2mm}
	\int\limits_{\partial T \backslash \Gamma} \hspace{-2mm}
	\vec{F}^{c,*}_n (\vec{U}_{h}^{\phantom{-}}\!,\, \vec{U}_{h}^{-}) \cdot \vec{V}_{h} \, d\bm{s}\\
	-
	\sum_{F \in \facets^{\text{int}}}
	\int\limits_{F}
	\avg{\mat{K}(\vec{U}_{h}) \nabla \vec{U}_{h}} : \jump{\vec{V}_{h}} \, d\bm{s}
	-
	\int\limits_{F}
	\avg{\vec{K}^T (\vec{U}_{h}) \nabla \vec{V}_{h}} : \jump{\vec{U}_{h}} \, d\bm{s}
	+
	\int\limits_{F}
	\uptau \avg{\vec{K} (\vec{U}_{h})} \jump{\vec{U}_{h}} : \jump{\vec{V}_{h}} \, d\bm{s},
\end{multline}
where $\vec{U}_{h}^{-}$ is the neighboring solution on element $T^-$ sharing the common facet  $F = \partial T \cap \partial T^- \in \facets^{\text{int}}$ and $\vec{F}^{c,*}_n$ denotes a \textit{numerical} flux, which stitches the discontinuous solutions into a uniquely-valued inviscid flux across the element boundaries.

In \Cref{eq::bilinear_form_interior_sdg}, $\uptau$ is the interior penalty coefficient that weakly enforces continuity, based on the viscous fluxes. Its expression for quadrilateral elements \cite{hillewaert2013development} is 
\begin{align} \label{eq::penalty_coefficient_ip}
	\uptau = \alpha (k + 1)^2 / l,
\end{align}
where $k$ is the polynomial order of the solution, $l$ is a characteristic length scale, taken as the ratio of the volume to the area of the boundary surface (of the element). In this work, we set $\alpha = 1$ for the DG, since a larger value requires a stricter time-step in an explicit time-marching scheme.

The terms in \Cref{eq::bilinear_form_interior_sdg} represent, from left to right, the volume contribution of the viscous and inviscid fluxes, the numerical flux term, the viscous surface term, the adjoint consistency term and the penalty term. Note, $\vec{K}$ is a \textit{block} matrix of viscous diffusion matrices, with viscous flux components given by $\vec{F}^{\nu}_i = \vec{K}_{ij} \nabla_j \vec{U}$, see \Cref{eq::definition_jacobians_inviscid_viscous,eq::diffusion_block_matrix_product}.

Throughout this work, the DG scheme employs the Vijayasundaram numerical flux \cite{vijayasundaram1986416}, defined as
\begin{align} \label{eq::upwinding_Vijayasundaram_sdg}
    \vec{F}^{c,*}_n ( \vec{U}^{+}_{h}, \vec{U}^{-}_{h} )
    =
    \mat{A}^{+}_{n} (\overline{\vec{U}}_{h})\, \vec{U}^{+}_{h} + \mat{A}^{-}_{n} (\overline{\vec{U}}_{h})\, \vec{U}^{-}_{h},
\end{align}
where upwinding is introduced through the eigenvalue decomposition of \Cref{subsec:inviscid_flux_splitting}, evaluated using the averaged state $\overline{\vec{U}}_{h} = \avg{\vec{U}_{h}}$.

\subsection{Hybridizable Discontinuous Galerkin: Interior Penalty}
\label{sec::hdg_ip}

Similar to \cite{franciolini2020104542}, we hybridize the DG method to derive a \textit{primal} HDG formulation. The resulting scheme is constructed in close analogy with the adjoint-symmetric DG-IP method presented in \Cref{sec::sdg}.

The resulting weak formulation of the HDG-IP scheme reads: find $\left( \vec{U}_{h}, \widehat{\vec{U}}_{h} \right) \in U_h \times \widehat{U}_h$, such that
\begin{subequations} \label{eq::weak_form_hdg_ip}
	\begin{align} \label{eq::weak_form_hdg_ip_element}
	\sum_{T \in \mesh} 
	\int\limits_{T} \pdt{\vec{U}_{h}} \cdot \vec{V}_{h} \, d\bm{x}
	+
	\mathcal{B}_{I} \big( \vec{U}_{h}, \widehat{\vec{U}}_{h}; \vec{V}_{h} \big)
	+
	\mathcal{B}_{\Gamma_{\!w}}
	&=
	0,
	\qquad 
	\forall\ \vec{V}_{h} \in U_{h},\\ \label{eq::weak_form_hdg_ip_facet}
	\widehat{\mathcal{B}}_{I} \big( \vec{U}_{h}, \widehat{\vec{U}}_{h}; \widehat{\vec{V}}_{h} \big)
	+
	\widehat{\mathcal{B}}_{\Gamma}
	&=
	0,
	\qquad 
	\forall\ \widehat{\vec{V}}_{h} \in \widehat{U}_{h},
\end{align}
\end{subequations}
where $(\cdot)_{I}$ and $(\cdot)_{\Gamma}$ denote interior terms and exterior/boundary terms, respectively, while the bilinear forms $\mathcal{B}$ and $\widehat{\mathcal{B}}$ correspond to functionals acting on the element test space $\vec{V}_{h}$ and facet test space $\widehat{\vec{V}}_{h}$, respectively.

In this section, the HDG-IP jump notation is defined as 
\begin{align}
	\widehat{\jump{\vec{U}_{h}}}
	:=  
	(\vec{U}_{h} - \widehat{\vec{U}}_{h}) \otimes \vec{n},
\end{align}
where $\vec{n}$ is the normal, taken with respect to the local element boundary $\partial T$.

In \Cref{eq::weak_form_hdg_ip}, the boundary forms $\mathcal{B}_{\Gamma}$, $\mathcal{B}_{\Gamma_{\!w}}$ and $\widehat{\mathcal{B}}_\Gamma$ are discussed in \ref{sec::hdg_ip_bcs}, while $\mathcal{B}_I$ is 
\begin{multline} \label{eq::bilinear_form_interior_hdg_ip_element}
	\mathcal{B}_{I} \big( \vec{U}_{h}, \widehat{\vec{U}}_{h}; \vec{V}_{h} \big)
	=
	\sum_{T \in \mesh} 
	\int\limits_{T}
	\Big( \vec{F}^{\nu}  (\vec{U}_{h}, \nabla \vec{U}_{h}) - \vec{F}^{c} (\vec{U}_{h}) \Big) : \grad{\vec{V}}_{h} \, d\bm{x}
	+
	\int\limits_{\partial T}
	\widehat{\vec{F}}^{c,*}_n (\vec{U}_{h}, \widehat{\vec{U}}_{h}) \cdot \vec{V}_{h} \, d\bm{s}\\
	-
	\hspace{-2mm}
	\int\limits_{\partial T \backslash \Gamma_{\!w}} \hspace{-2mm}
	\vec{K} (\widehat{\vec{U}}_{h}) \nabla \vec{U}_{h} : (\vec{V}_{h} \otimes \vec{n}) \, d\bm{s}
	-
	\hspace{-2mm}
	\int\limits_{\partial T \backslash \Gamma_{\!w}} \hspace{-2mm}
	\vec{K}^T (\widehat{\vec{U}}_{h}) \nabla \vec{V}_{h} : \widehat{\jump{\vec{U}_{h}}} \, d\bm{s}
	+
	\hspace{-2mm}
	\int\limits_{\partial T \backslash \Gamma_{\!w}} \hspace{-2mm}
	\uptau \vec{K} (\widehat{\vec{U}}_{h}) \widehat{\jump{\vec{U}_{h}}} : (\vec{V}_{h} \otimes \vec{n})\, d\bm{s},
\end{multline}
where the HDG-IP penalty coefficient $\uptau$ is defined using $\alpha = 10$ to ensure coercivity\footnote{Since this term is treated implicitly (we only consider implicit HDG), it does not impose any additional time-step restriction--unlike explicit time integration schemes.}, see \Cref{eq::penalty_coefficient_ip}.

The transmission conditions associated with the form $\widehat{\mathcal{B}}_I$ in \Cref{eq::weak_form_hdg_ip_facet} is defined as
\begin{align} \label{eq::bilinear_form_interior_hdg_ip_facet}
	\widehat{\mathcal{B}}_{I} \big( \vec{U}_{h}, \widehat{\vec{U}}_{h}; \widehat{\vec{V}}_{h} \big)
	=
	-
	\sum_{T \in \mesh} 
	\int\limits_{\partial T \backslash \Gamma} \hspace{-2mm}
	\widehat{\vec{F}}^{c,*}_n (\vec{U}_{h}, \widehat{\vec{U}}_{h}) \cdot \widehat{\vec{V}}_{h} \, d\bm{s}
	+
	\hspace{-2mm}
	\int\limits_{\partial T \backslash \Gamma}
	\hspace{-2mm}
	\Big( 
	\vec{K} (\widehat{\vec{U}}_{h}) \nabla \vec{U}_{h}
	-
	\uptau \vec{K} (\widehat{\vec{U}}_{h}) \widehat{\jump{\vec{U}_{h}}} 
	\Big) : ( \widehat{\vec{V}}_{h} \otimes \vec{n} )\, d\bm{s}.
\end{align}

In \Cref{eq::bilinear_form_interior_hdg_ip_element}, the terms from left to right represent the volume contribution of the viscous and inviscid fluxes, the numerical flux term, the viscous surface term, the adjoint consistent term and the penalty term. Similarly, in \Cref{eq::bilinear_form_interior_hdg_ip_facet}, the terms correspond to the numerical flux contribution, followed by the viscous surface term and its penalty term.

The HDG numerical (inviscid) flux is defined as
\begin{align} \label{eq::inviscid_numerical_flux}
	\widehat{\vec{F}}^{c,*}_n (\vec{U}_{h}, \widehat{\vec{U}}_{h})
	= 
	\vec{F}^{c}_{n} (\widehat{\vec{U}}_{h}) + \mat{S}_{c} (\widehat{\vec{U}}_{h}) (\vec{U}_{h} - \hat{\vec{U}}_{h}),
\end{align}
where we obtain an upwinded flux vector splitting strategy by specifying the convective stabilization matrix as $\mat{S}_{c} (\widehat{\vec{U}}_{h}) = \vec{A}_n^+ (\widehat{\vec{U}}_{h})$. By substituting the value of $\vec{S}_c$ and given that $\mat{A}_n = \mat{A}_{n}^{+} + \mat{A}_{n}^{-}$, with $\vec{F}^{c}_{n} = \mat{A}_{n} \vec{U}$, we get an analogous expression to the DG flux in \Cref{eq::upwinding_Vijayasundaram_sdg}, namely
\begin{align}
	\widehat{\vec{F}}^{c,*}_n (\vec{U}_{h}, \widehat{\vec{U}}_{h})
	=
	\mat{A}_{n}^{+} (\widehat{\vec{U}}_{h})\, \vec{U}_{h} + \mat{A}_{n}^{-} (\widehat{\vec{U}}_{h})\, \widehat{\vec{U}}_{h}.
\end{align}

\begin{remark}
	Compared to the DG-IP formulation in \Cref{eq::bilinear_form_interior_sdg}, the averages are replaced with the corresponding terms evaluated using the facet variables $\widehat{\vec{U}}$, since they are uniquely defined. For instance, $\avg{\mat{K}(\vec{U})}$ in \Cref{sec::sdg} is replaced with $\mat{K} (\widehat{\vec{U}})$ in \Cref{sec::hdg_ip}.
\end{remark}
\begin{remark}
	From \Cref{eq::weak_form_hdg_ip}, it is evident that the boundary conditions are enforced through $\widehat{\mathcal{B}}_\Gamma$ via the transmission conditions in \Cref{eq::weak_form_hdg_ip_facet}. The only exception is adiabatic walls, where we modify the viscous flux, see \Cref{eq::adiabatic_viscous_flux_compact}, by enforcing no heat flux in the normal direction in $\mathcal{B}_{\Gamma_{\!w}}$ within \Cref{eq::weak_form_hdg_ip_element}.
\end{remark}

\subsection{Hybridizable Discontinuous Galerkin: Mixed Formulation}
\label{sec::hdg_mixed}

An alternative HDG formulation \cite{vila-perezHybridisableDiscontinuousGalerkin2021} of the compressible Navier-Stokes equations, denoted by HDG-MX, is given by 
introducing additional variables as unknowns, namely the symmetric rate-of-strain tensor $\EPS$ and the temperature gradient $\vec{\phi}$, defined as
\begin{alignat}{2}
    \EPS & := \begin{pmatrix}
        \varepsilon_{11} & \varepsilon_{12} \\[0.8ex]
        \varepsilon_{12} & \varepsilon_{22}
    \end{pmatrix} = \frac{\grad{\VEL} + \grad{\VEL}^\T}{2} - \frac{1}{3} \nabla \cdot (\VEL) \I, & \qquad
     \vec{\phi} & := \begin{pmatrix}
        \phi_1 \\[0.8ex]
        \phi_2
    \end{pmatrix} = \grad{\theta}.
\end{alignat}

With these definitions at hand, we express the viscous fluxes from \Cref{eq::inviscid_viscous_fluxes} as a function of the conservative variables~$\vec{U}$, the rate-of-strain tensor $\EPS$ and the temperature gradient $\vec{\phi}$, as
\begin{alignat}{2}
	\vec{F}^{\nu}(\vec{U}, \EPS, \vec{\phi}) := \begin{pmatrix}
    0 & 0 \\[0.8ex]
    2 \mu \varepsilon_{11} & 2 \mu \varepsilon_{12} \\[0.8ex]
    2 \mu \varepsilon_{12} & 2 \mu \varepsilon_{22} \\[0.8ex]
    2 \mu (\varepsilon_{11} u_1 + \varepsilon_{12} u_2) + \kappa \phi_1 & 2 \mu (\varepsilon_{12} u_1 + \varepsilon_{22} u_2) + \kappa \phi_2    
\end{pmatrix}.
\end{alignat}

The weak formulation of the HDG-MX method reads as follows: Find $\left(\vec{U}_h,\hat{\vec{U}}_h, (\EPS_h, \vec{\phi}_h) \right) \in U_h \times \hat{U}_h \times Q_h$, such that

\begin{subequations}
\label{eq::weak_form_hdg_mixed}
\begin{align} \label{eq::weak_form_hdg_mixed_element}
	\sum_{T \in \mesh} 
	\int\limits_{T} \pdt{\vec{U}_h} \cdot \vec{V}_h \, d\bm{x}
	+
	\mathcal{B}_I \big( \vec{U}_h, \widehat{\vec{U}}_h, \EPS_h, \vec{\phi}_h; \vec{V}_h \big)
	&= 0,
	\qquad 
	\forall\ \vec{V} \in U_h,\\ \label{eq::weak_form_hdg_mixed_facet}
	\phantom{\sum_{T \in \mesh} \int\limits_{T}}
    \widehat{\mathcal{B}}_I \big( \vec{U}_h, \widehat{\vec{U}}_h, \EPS_h, \vec{\phi}_h; \widehat{\vec{V}}_h \big) +
	\widehat{\mathcal{B}}_{\Gamma}
	&= 0,
	\qquad 
	\forall\ \widehat{\vec{V}}_h \in \widehat{U}_h,\\ \label{eq::weak_form_hdg_mixed_auxiliary}
    \tilde{\mathcal{B}}_I^\varepsilon \big( \vec{U}_h, \widehat{\vec{U}}_h, \EPS_h; \mat{\zeta}_h \big) 
    + 
    \tilde{\mathcal{B}}_I^{\phi} \big( \vec{U}_h, \widehat{\vec{U}}_h, \vec{\phi}_h; \vec{\varphi}_h \big)  
    & = 0, 	\qquad \forall\ (\mat{\zeta}_h, \vec{\varphi}_h) \in Q_h,
\end{align}    
\end{subequations}
where the forms $\mathcal{B}_I$, $\widehat{\mathcal{B}}_I$, $\tilde{\mathcal{B}}^\varepsilon_I$ and $\tilde{\mathcal{B}}^{\phi}_I$ are defined as
\begin{subequations}
\label{eq::bilinear_form_hdg_mixed}
\begin{align}
    \begin{split}
        \mathcal{B}_I \big( \vec{U}_h, \widehat{\vec{U}}_h, \EPS_h, \vec{\phi}_h; \vec{V}_h \big) 
        &:=  
        \sum_{T \in \mesh}
        \int\limits_{T} \Big( \vec{F}^{\nu}  (\vec{U}_h, \EPS_h, \vec{\phi}_h) - \vec{F}^{c} (\vec{U}_h) \Big) : \grad{\vec{V}_h} \, d\bm{x} \\
         &+ 
         \int\limits_{\partial T} \left(\widehat{\vec{F}}^{c,*}_n (\vec{U}_h, \widehat{\vec{U}}_h) - \vec{F}_n^{\nu, *}  (\vec{U}_h, \widehat{\vec{U}}_h, \EPS_h, \vec{\phi}_h) \right) \cdot \vec{V}_h \, d\bm{s}, 
    \end{split}\\
    \widehat{\mathcal{B}}_I \big( \vec{U}_h, \widehat{\vec{U}}_h, \EPS_h, \vec{\phi}_h; \widehat{\vec{V}}_h \big) 
    &:= 
    \sum_{T \in \mesh} 
    \int\limits_{\partial T \backslash \Gamma} 
    \Big( \widehat{\vec{F}}^{c,*}_n (\vec{U}_h, \widehat{\vec{U}}_h) - \vec{F}_n^{\nu, *}  (\vec{U}_h, \widehat{\vec{U}}_h, \EPS_h, \vec{\phi}_h) \Big) \cdot \widehat{\vec{V}}_h \, d\bm{s}\\
    \begin{split}
    \tilde{\mathcal{B}}_I^\varepsilon \big( \vec{U}_h, \widehat{\vec{U}}_h, \EPS_h; \mat{\zeta}_h \big)
    &:= 
    \sum_{T  \in \mesh}
    \int\limits_{T}
    \EPS_h : \mat{\zeta}_h \, d\bm{x} 
    + 
    \int\limits_{T} \VEL_h \cdot \nabla \cdot \Big( \mat{\zeta}_h - \frac{1}{3}\tr(\mat{\zeta}_h)\I \Big) \, d\bm{x}\\ 
    &- 
    \int\limits_{\partial T} \hat{\VEL}_h \cdot \Big( \mat{\zeta}_h - \frac{1}{3}\tr(\mat{\zeta}_h)\I \Big) \vec{n} \, d\bm{s},   
    \end{split}\\
    \tilde{\mathcal{B}}_I^{\phi} \big( \vec{U}_h, \widehat{\vec{U}}_h, \vec{\phi}_h; \vec{\varphi}_h \big) 
    &:= 
    \sum_{T  \in \mesh} 
    \int\limits_{T} \vec{\phi}_h \cdot \vec{\varphi}_h \, d\bm{x} 
    + 
    \int\limits_{T} \theta_h \nabla \cdot (\vec{\varphi}_h) \, d\bm{x} 
    - 
    \int\limits_{\partial T} \hat{\theta}_h \vec{\varphi}_h \cdot \vec{n} \, d\bm{s}.
\end{align}    
\end{subequations}
Note that in \Cref{eq::bilinear_form_hdg_mixed}, the discrete velocities $\vec{u}_h := \vec{u}(\vec{U}_h)$, $\hat{\vec{u}}_h := \vec{u}(\hat{\vec{U}}_h)$, and the discrete temperatures $\theta_h := \theta(\vec{U}_h)$, $\hat{\theta}_h := \theta(\hat{\vec{U}}_h)$, are a function of the conservative variables $\vec{U}_h$ and $\hat{\vec{U}}_h$, respectively. Moreover, the inviscid numerical flux $\widehat{\vec{F}}^{c,*}_n$ is defined as in \Cref{eq::inviscid_numerical_flux}, whereas the viscous numerical flux $\vec{F}_n^{\nu, *}$ is given by
\begin{align}
    \label{eq::viscous_numerical_flux}
    \vec{F}_n^{\nu, *}  (\vec{U}_h, \widehat{\vec{U}}_h, \EPS_h, \vec{\phi}_h) & := \vec{F}_n^{\nu}(\hat{\vec{U}}_h, \EPS_h, \vec{\phi}_h) - \mat{S}_d (\vec{U}_h - \hat{\vec{U}}_h).
\end{align}
In this work we choose the viscous stabilization matrix $\mat{S}_d$ as $\diag(0, \mu, \mu, \kappa)$.

In a similar fashion as in the HDG-IP formulation of \Cref{sec::hdg_ip}, the boundary conditions are imposed in a weak sense by defining a suitable boundary operator $\widehat{\mathcal{B}}_{\Gamma}$ in \Cref{eq::weak_form_hdg_mixed}.The definition of
the boundary operator for an adiabatic wall or far-field boundary condition is given in \ref{sec::hdg_mixed_bcs}, while for more sophisticated boundary conditions, we refer to \cite{ellmenreichCharacteristicBoundaryConditions2026}.

\begin{remark}
    Unlike the HDG-IP formulation, HDG-MX attains coercivity without tuning of the penalty term. Additionally, by employing approximations of the same polynomial order $k$ for the primal, mixed, and hybrid unknowns yields optimal convergence rates for the viscous stress and heat flux \cite{vila-perezHybridisableDiscontinuousGalerkin2021}.
\end{remark}

\section{Implicit-Explicit Time-Marching Strategy}
\label{sec::imex}

In this work, we consider a geometry-split IMEX strategy in which the mesh $\mesh$ is partitioned into two submeshes $\mesh^{im} \subseteq \mesh$ and $\mesh^{ex} \subseteq \mesh$, such that
\begin{align}
	\label{eq::union_mesh}
	\mesh &= \mesh^{im} \cup \mesh^{ex},
\end{align}
with interface $\Gamma_{\!i} = \mesh^{im} \cap \mesh^{ex}$. Let $n$ denote the total number of elements in $\mesh$ and let $n_{i}$ and $n_{e}$ denote the number of implicit and explicit elements in $\mesh^{im}$ and $\mesh^{ex}$, respectively. Then, from \Cref{eq::union_mesh} follows the relationship between the total number of elements as 
\begin{align}
	\label{eq::number_elements}
	n = n_{i} + n_{e}.
\end{align}

Unless stated otherwise, on mesh $\mesh^{im}$ we discretize in space using an HDG method from \Cref{sec::hdg_ip,sec::hdg_mixed} and treat the resulting system of equations implicitly in time, while on $\mesh^{ex}$ we discretize in space using a DG method from \Cref{sec::sdg} and treat the resulting system of equations explicitly in time. We denote by $(\cdot)^{im}$ and $(\cdot)^{ex}$ the discrete quantities defined on the partitions $\mesh^{im}$ and $\mesh^{ex}$, respectively.



\subsection{Interface Conditions}
\label{subsec::interface_conditions}
The solutions on $\mesh^{im}$ and $\mesh^{ex}$ are weakly coupled across the interface $\Gamma_{\!i}$ through the imposition of suitable interface conditions. Furthermore, both solutions are advanced in a time-synchronized manner, which ensures the conservation of the coupled system.



\subsubsection{Coupling from the perspective of a DG scheme}
\label{subsec::coupling_sdg}

In the DG-IP formulation, the interface is treated in the same manner as any interior face. Consequently, on the interface $\Gamma_{\!i}$, the boundary contribution appearing in \Cref{eq::weak_form_sdg_ip} takes the form
\begin{multline} \label{eq::interface_bilinearform_explicit}
	\mathcal{B}_{\Gamma} 
	:= 
	\mathcal{B}_{\Gamma_{\!i}} \big( \vec{U}_{h}, \vec{U}^{im}_{h}; \vec{V}_{h} \big)
	=
	\sum_{F \in \Gamma_{\!i}}
	\int\limits_{F}
	\vec{F}^{c,*}_n (\vec{U}_{h}, \vec{U}^{im}_{h}) \cdot \vec{V}_{h} \, d\bm{s}
	-
	\int\limits_{F}
	\avg{\vec{K} (\vec{U}_{h}) \nabla \vec{U}_{h} }_{\Gamma_{\!i}} : (\vec{V}_{h} \otimes \vec{n}) \, d\bm{s}\\
	-
	\int\limits_{F}
	\avg{\vec{K}^T (\vec{U}_{h})}_{\Gamma_{\!i}} \nabla \vec{V}_{h} : \jump{\vec{U}_{h}}_{\Gamma_{\!i}} \, d\bm{s}
	+
	\int\limits_{F}
	\uptau \avg{\vec{K} (\vec{U}_{h})}_{\Gamma_{\!i}} \jump{\vec{U}_{h}}_{\Gamma_{\!i}} : (\vec{V}_{h} \otimes \vec{n}) \, d\bm{s},
\end{multline}
with interface average and jump operators $\avg{\vec{\psi}}_{\Gamma_{\!i}} = \left(\vec{\psi} + \vec{\psi}^{im}\right)/2$ and $\jump{\vec{\psi}}_{\Gamma_{\!i}} = (\vec{\psi} - \vec{\psi}^{im}) \otimes \vec{n}$.

\subsubsection{Coupling from the perspective of an HDG scheme}
\label{subsec::coupling_hdg}

For both HDG-IP and HDG-MX formulations, the interface condition is imposed via transmissibility conditions expressed as a Dirichlet constraint. While this condition is enforced strongly on the facet variables, it induces only a weak imposition on the element-local solution $\vec{U}_h$. Hence, on an interface $\Gamma_{\!i}$, the boundary contribution in \Cref{eq::weak_form_hdg_ip_facet,eq::weak_form_hdg_mixed_facet} is given by
\begin{align}  \label{eq::interface_bilinearform_implicit_facet}
	\widehat{\mathcal{B}}_{\Gamma} 
	:= 
	\widehat{\mathcal{B}}_{\Gamma_{\!i}} \big( \widehat{\vec{U}}_{h}, \vec{U}^{ex}_{h}; \widehat{\vec{V}}_{h} \big)
	=
	\sum_{F \in \Gamma_{\!i}}
	\int\limits_{F}
	(\widehat{\vec{U}}_{h} - \vec{U}^{ex}_{h}) \cdot \widehat{\vec{V}}_{h} \, d\bm{s}. 
\end{align}

\begin{remark}
	The interface coupling in \Cref{eq::interface_bilinearform_implicit_facet} for HDG formulations depends solely on the external conservative variables, regardless of the specific HDG variant used (e.g., mixed, interior penalty or inviscid).
\end{remark}
	
\begin{remark}
	For the DG formulation, the coupling in \Cref{eq::interface_bilinearform_explicit} depends on the type of external discretization. In particular, when coupling with the mixed formulation in \Cref{sec::hdg_mixed}, the viscous flux involves both conservative and auxiliary variables and is evaluated as
	\begin{align}
		\avg{\vec{K} (\vec{U}_{h}) \nabla \vec{U}_{h} }_{\Gamma_{\!i}}
		:= \frac{1}{2} 
		\Big(
		\vec{F}^{\nu} (\vec{U}_h, \nabla \vec{U}_h)
		+ 
		\vec{F}^{\nu} (\vec{U}_{h}^{im},\, \EPS_h^{im},\, \vec{\phi}_h^{im})
		\Big).
	\end{align}
\end{remark}

\subsection{Temporal Discretization}
\label{subsec::temporal_discretization}

\subsubsection{Semi-discrete form}

After spatial discretization, the governing PDE reduces to a system of ordinary differential equations. The explicit (DG) and implicit (HDG) parts of the semi-discrete system read
\begin{subequations} \label{eq::semidiscrete_time}
	\begin{align} \label{eq::semidiscrete_time_sdg}
		\mat{M} \frac{\partial \vec{U}^{ex}\hspace{-3mm}}{\partial t} \ - \
		\vec{R}^{\mathrm{DG}} \left( \vec{U}^{ex},\, \vec{U}^{im} \right) &= \vec{0},\\ \label{eq::semidiscrete_time_hdg_1}
		\mat{M} \frac{\partial \vec{U}^{im}\hspace{-3mm}}{\partial t} \ - \
		\vec{R}^{\mathrm{HDG}} \left( \vec{U}^{im},\, \widehat{\vec{U}}^{im},\, \vec{Q}^{im} \right) &= \vec{0},\\ \label{eq::semidiscrete_time_hdg_2}
		\phantom{\mat{M} \frac{\partial \vec{U}^{ex}\hspace{-3mm}}{\partial t}}
		\widehat{\vec{R}}^{\mathrm{HDG}} \left( \vec{U}^{im},\, \widehat{\vec{U}}^{im},\, \vec{Q}^{im},\, \vec{U}^{ex} \right) &= \vec{0},\\ \label{eq::semidiscrete_time_hdg_3}
		\phantom{\mat{M} \frac{\partial \vec{U}^{ex}\hspace{-3mm}}{\partial t}}
		\widetilde{\vec{R}}^{\mathrm{HDG}} \left( \vec{U}^{im},\, \widehat{\vec{U}}^{im},\, \vec{Q}^{im} \right) &= \vec{0},
	\end{align}
\end{subequations}
where $\mat{M}$ denotes the mass matrix and the HDG mixed variables\footnote{In the case of an HDG-IP, no mixed variables are needed, i.e. $\vec{Q} = \vec{0}$.} are compactly defined as $\vec{Q} = \left( \EPS,\, \vec{\phi} \right)$.

In \Cref{eq::semidiscrete_time}, $\vec{R}$ denotes the physical residual vector arising from the spatial discretization. It is defined through the bilinear forms in \Cref{eq::weak_form_sdg_ip,eq::weak_form_hdg_ip,eq::weak_form_hdg_mixed} and incorporates the coupling between the implicit and explicit solutions across the interface $\Gamma_{\!i}$. For brevity, the dependence of $\vec{R}$ on the physical boundary states $\vec{U}_{\infty}$ and $\vec{U}_{w}$ is omitted. The interface coupling terms in \Cref{eq::semidiscrete_time_sdg,eq::semidiscrete_time_hdg_2} are obtained from \Cref{eq::interface_bilinearform_explicit,eq::interface_bilinearform_implicit_facet}, respectively.

\subsubsection{Multi-stage Additive Runge-Kutta Methods}

A generic Butcher tableau for a pair of additive Runge-Kutta (ARK) methods is shown in \Cref{tab::butcher_tableau_sidebyside}. In this framework, the implicit part of the scheme is constructed using a singly diagonally implicit Runge-Kutta (SDIRK) method, while the explicit part is represented by a standard explicit Runge-Kutta (ERK) method.
\begin{table}[tb!] 
		\begin{minipage}{\textwidth}
			\centering
	\begin{tabular}{c!{\vrule width 1pt}ccccc}
			\multicolumn{6}{c}{\,} \\
			0                    & 0                      & 0                      & $\cdots$ & 0                      & 0 \\ 
			$\overline{c}_2$     & $\overline{a}_{21}$    & 0                      & $\cdots$ & 0                      & 0 \\
			$\overline{c}_3$     & $\overline{a}_{31}$    & $\overline{a}_{32}$    & $\cdots$ & 0                      & 0 \\
			$\vdots$             & $\vdots$               & $\vdots$               & $\ddots$ & $\vdots$               & $\vdots$ \\
			$\overline{c}_{s+1}$ & $\overline{a}_{s+1,1}$ & $\overline{a}_{s+1,2}$ & $\cdots$ & $\overline{a}_{s+1,s}$ & 0 \\ \noalign{\hrule height 1pt}
			\textbf{ERK}         & $\overline{b}_1$       & $\overline{b}_2$       & $\cdots$ & $\overline{b}_s$       & $\overline{b}_{s+1}$
		\end{tabular}
		\hspace{2em} 
		\begin{tabular}{c!{\vrule width 1pt}ccccc}
			\multicolumn{6}{c}{\,} \\ 
			0             & 0        & 0        & 0        & $\cdots$ & 0 \\
			$c_1$         & 0        & $a_{11}$ & 0        & $\cdots$ & 0 \\
			$c_2$         & 0        & $a_{21}$ & $a_{22}$ & $\cdots$ & 0 \\
			$\vdots$      & $\vdots$ & $\vdots$ & $\ddots$ & $\ddots$ & $\vdots$ \\
			$c_s$         & 0        & $a_{s1}$ & $a_{s2}$ & $\cdots$ & $a_{ss}$ \\ \noalign{\hrule height 1pt}
			\textbf{DIRK} & 0        & $b_1$    & $b_2$    & $\cdots$ & $b_s$
		\end{tabular}
	\end{minipage}
	\caption{ARK methods: Explicit (ERK) and Implicit (DIRK) Butcher tableaux synchronized via padding.}
	\label{tab::butcher_tableau_sidebyside}
\end{table}

The current IMEX scheme is synchronized via padding, meaning that the stage times of the explicit and implicit methods coincide, i.e. $\overline{c}_{i+1} = c_i$. This choice follows the classical structure of an ARK method, which ensures that both parts of the scheme advance in time consistently and facilitates the coupling between explicit and implicit components. For the numerical experiments in this work, we adopt ARS-type schemes~\cite{ascher1997implicit}, which are well-established IMEX methods. To update the solution from time $t_{n}$ to $t_{n+1}$, refer to \Cref{alg::imex_scheme_update}. 

The explicit domain is advanced using an $(s+1)$-stage explicit Runge-Kutta scheme
\begin{subequations}
	\begin{align} \label{eq::erk_solution_update}
		\vec{U}_{n+1}^{ex} 
		&=
		\vec{U}_{n}^{ex}
		-
		\Delta t\, \mat{M}^{-1} \sum_{i=1}^{s+1} \overline{b}_{i}\,
		\vec{R}^{\mathrm{DG}} \left( \vec{y}_{i}^{ex},\, \vec{y}_{i-1}^{im} \right),\\ \label{eq::erk_stage_update}
		\vec{y}_{i}^{ex}
		&=
		\begin{cases}
			\vec{U}_{n}^{ex}, & i = 1,\\
			\vec{U}_{n}^{ex}
			-
			\Delta t\, \mat{M}^{-1} \displaystyle\sum_{j=1}^{i-1} \overline{a}_{ij}\,
			\vec{R}^{\mathrm{DG}} \left( \vec{y}_{j}^{ex},\, \vec{y}_{j-1}^{im} \right), & i = 2,\, \dots,\, s\!+\!1,
		\end{cases}
\end{align}
\end{subequations}
where $\Delta t$ is the time step, $\vec{y}^{ex}$ and $\vec{y}^{im}$ are the intermediate stage solutions for the explicit and implicit schemes, respectively. Note, stage zero of the implicit solution is shorthand for $\vec{y}_{0}^{im} = \vec{U}_{n}^{im}$.

The implicit domain is advanced using an $s$-stage SDIRK scheme
\begin{subequations}
	\begin{align} \label{eq::dirk_solution_update}
		\vec{U}_{n+1}^{im}
		&=
		\vec{U}_{n}^{im}
		-
		\Delta t\, \mat{M}^{-1} \sum_{i=1}^{s} b_{i}\,
		\vec{R}^{\mathrm{HDG}} \left( \vec{y}_{i}^{im},\, \widehat{\vec{y}}_{i}^{im} \right),
		\\ \label{eq::dirk_stage_update_solution}
		\mat{M} \vec{y}_{i}^{im}
		+
		\Delta t\, a_{ii}\, 
		\vec{R}^{\mathrm{HDG}} \left( \vec{y}_{i}^{im},\, \widehat{\vec{y}}_{i}^{im},\, \widetilde{\vec{y}}_{i}^{im} \right)
		&=
		\mat{M} \vec{U}_{n}^{im}
		-
		\Delta t\, \sum_{j=1}^{i-1} a_{ij}\, 
		\vec{R}^{\mathrm{HDG}} \left( \vec{y}_{j}^{im},\, \widehat{\vec{y}}_{j}^{im},\, \widetilde{\vec{y}}_{j}^{im} \right),
		\\ \label{eq::dirk_stage_update_facet}
		\widehat{\vec{R}}^{\mathrm{HDG}} \left( \vec{y}_{i}^{im},\, \widehat{\vec{y}}_{i}^{im},\, \widetilde{\vec{y}}_{i}^{im},\, \vec{y}_{i+1}^{ex} \right)
		&=
		\vec{0},
		\\[8pt] \label{eq::dirk_stage_update_mixed}
		\widetilde{\vec{R}}^{\mathrm{HDG}} \left( \vec{y}_{i}^{im},\, \widehat{\vec{y}}_{i}^{im},\, \widetilde{\vec{y}}_{i}^{im} \right)
		&=
		\vec{0},
	\end{align}
\end{subequations}
where $\widehat{\vec{y}}$ and $\widetilde{\vec{y}}$ correspond to the stage solutions for the facet $\widehat{\vec{U}}$ and mixed variables $\vec{Q}$, respectively.

\subsection{Global Implicit System}
\label{subsec::global_implicit_system}

The global system in \Cref{eq::dirk_stage_update_solution,eq::dirk_stage_update_facet,eq::dirk_stage_update_mixed} is solved using a Newton-Raphson iterative procedure. Let $m$ denote the nonlinear iteration index. At each iteration, we seek updates to the solution variables in the form of increments
\begin{align}
	\delta \vec{U}^{m} = \vec{U}^{m+1} - \vec{U}^{m},
	\qquad
	\delta \vec{Q}^{m} = \vec{Q}^{m+1} - \vec{Q}^{m},
	\qquad
	\delta \vec{\widehat{U}}^{m} = \vec{\widehat{U}}^{m+1} - \vec{\widehat{U}}^{m}.
\end{align}

At iteration $m$, the nonlinear system is linearized about the current iterate, yielding a Jacobian system for the increments. Denoting by $\vec{f}^{m}$, $\vec{g}^{m}$ and $\vec{h}^{m}$ the residuals associated with \Cref{eq::dirk_stage_update_solution,eq::dirk_stage_update_facet,eq::dirk_stage_update_mixed}, respectively, the linearized system can be written in block form as
\begin{align} \setstretch{1.3} \label{eq::implicit_system}
	\begin{pmatrix}
		\begin{array}{cc|c}
			\mathbb{A}_{11}^{m} & \mathbb{A}_{12}^{m} & \mathbb{A}_{13}^{m}\\
			\mathbb{A}_{21}^{m} & \mathbb{A}_{22}^{m} & \mathbb{A}_{23}^{m}\\ \hline
			\mathbb{A}_{31}^{m} & \mathbb{A}_{32}^{m} & \mathbb{A}_{33}^{m}
		\end{array}
	\end{pmatrix}
	\begin{pmatrix}
		\begin{array}{c}
			\delta \vec{U}^{m}\\
			\delta \vec{Q}^{m}\\ \hline
			\delta \vec{\widehat{U}}^{m}
		\end{array}
	\end{pmatrix}
	=
	\begin{pmatrix}
		\begin{array}{c}
			\vec{f}^{m}\\
			\vec{g}^{m}\\ \hline
			\vec{h}^{m}
		\end{array}
	\end{pmatrix}
	\qquad
	\Longleftrightarrow
	\qquad
	\begin{pmatrix}
		\mathbb{K}^{m} & \mathbb{C}^{m}\\
		\mathbb{D}^{m} & \mathbb{A}_{33}^{m}
	\end{pmatrix}
	\begin{pmatrix}
		\delta \vec{Z}^{m}\\
		\delta \vec{\widehat{U}}^{m}
	\end{pmatrix}
	=
	\begin{pmatrix}
		\vec{w}^{m}\\
		\vec{h}^{m}
	\end{pmatrix}.
\end{align}

For convenience, the compact block notation in \Cref{eq::implicit_system} is defined as
\begin{align} \setstretch{1.3}
	\mathbb{K}^{m}
	=
	\begin{pmatrix}
		\mathbb{A}_{11}^{m} & \mathbb{A}_{12}^{m}\\
		\mathbb{A}_{21}^{m} & \mathbb{A}_{22}^{m}
	\end{pmatrix},
	\quad
	\mathbb{C}^{m}
	=
	\begin{pmatrix}
		\mathbb{A}_{13}^{m}\\
		\mathbb{A}_{23}^{m}
	\end{pmatrix},
	\quad
	\mathbb{D}^{m}
	=
	\begin{pmatrix}
		\mathbb{A}_{31}^{m} & \mathbb{A}_{32}^{m}
	\end{pmatrix},
	\quad
	\delta \vec{Z}^{m}
	=
	\begin{pmatrix}
		\delta \vec{U}^{m}\\
		\delta \vec{Q}^{m}
	\end{pmatrix},
	\quad
	\vec{w}^{m}
	=
	\begin{pmatrix}
		\vec{f}^{m}\\
		\vec{g}^{m}
	\end{pmatrix}.
\end{align}

The block structure highlights the separation between element-interior unknowns, $\delta \vec{Z}^{m}$, and globally coupled trace unknowns, $\delta \vec{\widehat{U}}^{m}$. In particular, the matrix $\mathbb{K}^{m}$ corresponds to purely local (element-wise) interactions, while $\mathbb{A}_{33}^{m}$ encodes the coupling between facet degrees of freedom across the mesh. The off-diagonal blocks $\mathbb{C}^{m}$ and $\mathbb{D}^{m}$ represent the interaction between interior and trace variables. Furthermore, since $\mathbb{K}^{m}$ is block-diagonal with respect to the mesh elements, it can be inverted locally and efficiently. This property enables the static condensation of the interior unknowns, leading to a reduced system posed solely in terms of the trace variable.

The above system can thus be statically condensed to obtain
\begin{align} \label{eq::statically_condensed_system}
	\mathbb{S}^{m} \delta \vec{\widehat{U}}^{m} = \vec{h}^{m} - \mathbb{D}^{m} \left[ \mathbb{K}^{m} \right]^{-1} \vec{w}^{m},
\end{align}
where the Schur complement of $\mathbb{K}$ is defined as
\begin{align} \label{eq::SchurCompliment}
	\mathbb{S}^{m} := \mathbb{A}_{33}^{m} - \mathbb{D}^{m} \left[ \mathbb{K}^{m} \right]^{-1} \mathbb{C}^{m}.
\end{align}

This condensed system is significantly smaller than the original global system, as it only involves degrees of freedom associated with the mesh skeleton. Once $\delta \vec{\widehat{U}}^{m}$ is obtained, the interior updates $\delta \vec{Z}^{m}$ can be recovered locally, element-by-element, via
\begin{align} \label{eq::local_solve_static_condensation}
	\delta \vec{Z}^{m}
	=
	\left[ \mathbb{K}^{m} \right]^{-1}
	\left(
	\vec{w}^{m}
	-
	\mathbb{C}^{m} \delta \vec{\widehat{U}}^{m}
	\right).
\end{align}

This two-step procedure--comprising a global solve for the trace unknowns followed by a local reconstruction of the interior fields--is a key feature of hybridized methods \cite{cockburn2009unified,cockburnStaticCondensationHybridization2016}. By reducing the size of the globally coupled system in \Cref{eq::statically_condensed_system}, it lowers the overall solution time and facilitates efficient parallelization of the local solves in \Cref{eq::local_solve_static_condensation}.

\begin{remark} \label{rem::hdg_cost}
	While the dimensions of the Schur complement and the facet unknowns are identical in HDG-IP and HDG-MX, the computational cost of both the Schur complement assembly in \Cref{eq::SchurCompliment} and the element-local reconstruction in \Cref{eq::local_solve_static_condensation} scales with the number of auxiliary variables. In HDG-IP, no mixed variables are introduced, leading to smaller local systems and reduced cost. In contrast, HDG-MX includes additional mixed variables, increasing the size of the local problems and thus the cost of both assembly and local solves.
\end{remark}

\subsection{Computational Performance}
\label{subsec::computational_efficiency}

Let $\mathcal{C}_{\textrm{total}}$ denote the total wall-clock time required by a (multi-stage) time integration scheme to advance the solution of polynomial degree $k$ over $(0, t_{\textrm{total}}]$ on a mesh $\mesh$, see \Cref{eq::union_mesh}, with $n$ elements. Assuming an average cost per time step $\mathcal{C}$, the total cost reads
\begin{align}
	\label{eq::wall_clock_time}
	\mathcal{C}_{\textrm{total}} = \frac{t_{\textrm{total}}}{\Delta t_s(\mesh, \vec{U}_h)} \, \mathcal{C}.
\end{align}

Here, $\Delta t_s(\mesh, \vec{U}_h)$ is the maximum stable time step of the (multi-stage) explicit scheme. It depends on both the mesh $\mesh$ and the local solution state $\vec{U}_h$, which encapsulates all relevant physical properties such as the speed of sound and viscosity. The dependence on the polynomial order is omitted for clarity; see \cite{hesthaven2008nodal,chalmers2020109095} for details.

For explicit schemes, the cost per time step scales linearly with the number of elements $n$,
\begin{align}
	\label{eq::explicit_wall_clock_time}
	\mathcal{C}^{ex}(n) \approx n \, \overline{\mathcal{C}}^{ex},
\end{align}
where $\overline{\mathcal{C}}^{ex}$ is the average cost per element.

For implicit schemes, the cost per time step depends on both the number of elements $n$ and the number of nonlinear iterations $m$. The dominant contributions are the assembly and the solution of the linearized system, yielding
\begin{align}
	\label{eq::implicit_wall_clock_time}
	\mathcal{C}^{im}(n, m) \approx m \left( \overline{\mathcal{C}}^{im}_{\textrm{assem}}(n) + \overline{\mathcal{C}}^{im}_{\textrm{solve}}(n) \right).
\end{align}

In \Cref{eq::implicit_wall_clock_time}, the precise scaling with $n$ depends on the solver, sparsity pattern and discretization. In the present setting, HDG discretizations are expected to reduce $\mathcal{C}^{im}(n,m)$ compared to standard DG methods due to the smaller globally coupled system resulting from static condensation (see \Cref{subsec::global_implicit_system}).

For an IMEX scheme, the cost per time step is decomposed into implicit and explicit contributions:
\begin{align}
	\label{eq::imex_wall_clock_time}
	\mathcal{C}^{imex}(n_i, n_e, m) := \mathcal{C}^{im}(n_i, m) + n_e \, \overline{\mathcal{C}}^{ex},
\end{align}
where $n_i + n_e = n$.

We define the speedup $S$ as the ratio of the total wall-clock time of a fully explicit scheme to that of the IMEX scheme:
\begin{align}
	\label{eq::speedup}
	S := \frac{\mathcal{C}^{ex}_{\textrm{total}}}{\mathcal{C}^{imex}_{\textrm{total}}}
	\approx 
	\frac{\Delta t_s^{imex}(\mesh^{ex}, \vec{U}_h)}{\Delta t_s^{ex}(\mesh, \vec{U}_h)}
	\cdot
	\frac{(n_i + n_e)\,\overline{\mathcal{C}}^{ex}}{\mathcal{C}^{im}(n_i,m) + n_e \overline{\mathcal{C}}^{ex}}
	= 
	r_t \cdot \eta,
\end{align}
where $r_t$ is the ratio of stable time steps between the IMEX and fully explicit scheme,
\begin{align}
	\label{eq::time_step_ratio}
	r_t := \frac{\Delta t_s^{imex}(\mesh^{ex}, \vec{U}_h)}{\Delta t_s^{ex}(\mesh, \vec{U}_h)},
\end{align}
and $\eta$ is the mean relative efficiency of the IMEX scheme,
\begin{align}
	\label{eq::efficiency}
	\eta := \frac{(n_i + n_e)\,\overline{\mathcal{C}}^{ex}}{\mathcal{C}^{im}(n_i, m) + n_e \overline{\mathcal{C}}^{ex}}.
\end{align}

Note that the stable time step of an IMEX scheme, $ \Delta t_s^{imex}(\mesh^{ex}, \vec{U}_h) $, depends on the explicit mesh $\mesh^{ex}$, rather than the entire mesh $\mesh$.

In the limiting cases of a fully implicit or fully explicit discretization, we obtain
\begin{equation}
	r_t =
	\begin{cases}
		\Delta t^{im}_{c} / \Delta t_s^{ex}(\mesh, \vec{U}_h),
		& \text{as } \mesh^{im} \to \mesh, \\
		1,
		& \text{as } \mesh^{ex} \to \mesh,
	\end{cases}
	\qquad
	\eta =
	\begin{cases}
		n_i \overline{\mathcal{C}}^{ex} \! /\, \mathcal{C}^{im}(n_i, m),
		& \text{as } \mesh^{im} \to \mesh, \\
		1,
		& \text{as } \mesh^{ex} \to \mesh,
	\end{cases}
\end{equation}
where $\Delta t^{im}_{c}$ is a \textit{cutoff} time step determined by accuracy requirements or physical time scales of interest.

The main purpose of adopting an IMEX strategy is to increase the stable time step, which is quantified by the ratio of stable time steps $r_t$. By construction $r_t \geq 1$, since
\begin{align}
\Delta t_s^{ex}(\mesh, \vec{U}_h) &\leq \Delta t_s^{imex}(\mesh^{ex}, \vec{U}_h),  \qquad \forall\ \mesh^{ex} \subseteq \mesh.
\end{align}

This gain is offset by the additional cost of the implicit solve, captured by the mean relative efficiency~$\eta$. As implicit methods are more expensive per element, we expect $n_i \overline{\mathcal{C}}^{ex\!} \! < \mathcal{C}^{im}(n_i, m)$, and therefore $\eta < 1$. Hence, $\eta$ quantifies the mean relative efficiency loss of the IMEX scheme compared to a fully explicit method.

From \Cref{eq::speedup}, a speedup $S > 1$ is achieved if the increase in stable time step outweighs the additional implicit cost, namely 
$r_t > 1 / \eta$. This condition highlights the importance of two factors in viable IMEX schemes: an \textit{effective mesh partitioning} and an \textit{efficient implicit solver}, to increase $r_t$ and $\eta$, respectively.

\begin{remark}
	A common strategy to reduce the implicit cost is to freeze the Jacobian matrix across iterations or time steps. While this lowers assembly and factorization costs, it may degrade robustness for strongly nonlinear problems or large time steps. This approach is not considered here in order to isolate the fundamental behavior of the IMEX formulation.
\end{remark}

\section{Numerical experiments}
\label{sec::experiments}

The IMEX splitting proposed in \Cref{sec::imex} has been assessed in a series of numerical experiments
for the compressible Navier-Stokes equations, including verification, validation and computational performance assessment. 
Our implementation is available in the open-source Python extension \texttt{DreAm} based on the finite element library \texttt{NGSolve} \cite{schoberlC++11ImplementationFinite2014} 
and the results of these experiments are presented in the following sections.

For both the HDG and DG spatial discretizations, the polynomial order is set to $k=3$, unless specified otherwise. We perform
an over-integration of the nonlinear terms using a quadrature rule up to polynomial order $3k$.

We denote by IMEX-IP the IMEX scheme based on the HDG-IP formulation, and by IMEX-MX the one based on the HDG-MX formulation. If appropriate, we 
include the results of a fully explicit DG-IP scheme and a fully implicit HDG scheme as a reference, since they use the same spatial discretization and time integration schemes as their IMEX counterparts. 

The nondimensionalization of Equations~\eqref{eq::compressible_equations} is carried out by using aerodynamic reference values for the velocity $u_r$, length $l_r$, density $\rho_r$, viscosity $\mu_r$ and temperature $T_r$.

\subsection{Verification and Validation}

\subsubsection{Method of Manufactured Solutions}
\label{sec::experiments::manufactured_solutions}
\pgfplotsset{table/search path={./figures/mms_data/}}

To verify the correctness of our IMEX implementation, we employ the method of manufactured
solutions (MMS) \cite{royVerificationEulerNavier2004,gassnerDiscontinuousGalerkinScheme2008}. This approach involves the construction of
a solution $\vec{U}_e$, hence the name manufactured solution, which is not necessarily a solution to the original PDEs. To impose the manufactured solution,
the governing equations are modified by introducing source terms $\vec{S}(\vec{U}_e, \nabla \vec{U}_e)$ and appropriate boundary conditions, i.e.,
\begin{alignat}{2}
    \pdt{\vec{U}} + \nabla \cdot ( \vec{F}^{c} - \vec{F}^{\nu}) & = \vec{S}(\vec{U}_e, \nabla \vec{U}_e) \quad \text{in} \quad \Omega \times (0, t_{\textrm{total}}] , \\
    \vec{U}                                                                   & = \vec{U}_e \quad \text{on} \quad \Gamma \times (0, t_{\textrm{total}}] ,                            \\
    \vec{S}(\vec{U}_e, \nabla \vec{U}_e)                                      & := \pdt{\vec{U}_e} + \nabla \cdot (\vec{F}(\vec{U}_e) - \vec{G}(\vec{U}_e, \nabla \vec{U}_e)),
\end{alignat}
that ensure the discrete solution $\vec{U}$ to satisfies the modified equations.

In this work, we consider the manufactured solution presented in \cite{gassnerDiscontinuousGalerkinScheme2008}, defined as
\begin{alignat}{2}
    \vec{U}_e & := \begin{pmatrix}
                       \sin(\pi (k_x x + k_y y) - \omega \pi t) + 4 \\
                       \sin(\pi (k_x x + k_y y) - \omega \pi t) + 4 \\
                       \sin(\pi (k_x x + k_y y) - \omega \pi t) + 4 \\
                       (\sin( \pi (k_x x + k_y y) - \omega \pi t) + 4)^2
                   \end{pmatrix},
\end{alignat}
over the domain $\Omega = [0,1]^2$ and time interval $\Delta T = (0, 0.01]$. The mesh $\mesh$ consists of $32 \times 32$ uniform quadrilateral elements,
with periodic boundary conditions applied on the boundary $\Gamma := \partial \Omega$ and the interface boundary conditions applied on $\Gamma_{\!i}$. We set $k_x = k_y = 2$ and $\omega=200$ to guarantee a periodic manufactured
solution in both space and time. The high temporal frequency $\omega$ and the fine mesh size have been chosen to ensure a dominant temporal discretization error over the spatial discretization error. The local Mach number is approximately $\Ma \in [0.95, 1.34]$, while setting the dynamic viscosity to $\mu = 0.01$ as in \cite{gassnerDiscontinuousGalerkinScheme2008} results in a local Reynolds number of approximately $\Re \in [424.26, 707.11]$.

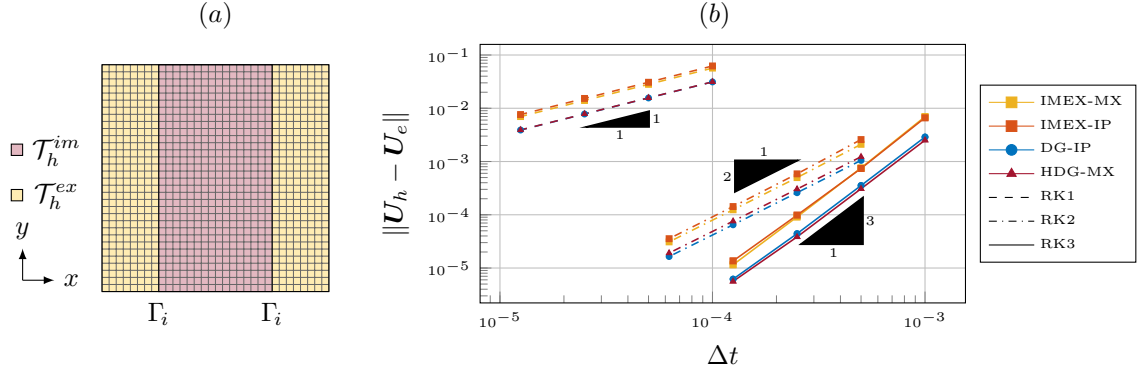
\begin{figure}[ht]
    \centering
    \tikzsetnextfilename{domain_mms}
\tikzsetnextfilename{domain_mms}
\begin{tikzpicture}[>=latex]
    \begin{scope}[scale=3, yshift=0.05cm]
          \fill[amber!30]  (0.0,0) rectangle (0.25,1);
    \fill[matred!30]  (0.25,0) rectangle (0.75,1);
    \fill[amber!30]  (0.75,0) rectangle (1.0,1);
    \draw[black!60, ultra thin] (0,0) grid[step=0.03125] (1,1);
    \draw[black, thin] (0,0) rectangle (1,1);
    \draw[black, thin] (0.25,0) rectangle (0.75,1);

    \begin{scope}[shift={(-0.4, 0.4)}]
        \fill[amber!30]  (0.0,0.0) rectangle (0.05, 0.05);
        \draw[black, thin] (0.0,0.0) rectangle (0.05, 0.05);
        \node at (0.2,0.02)  {$\mesh^{ex}$};
    \end{scope}

    \begin{scope}[shift={(-0.4, 0.6)}]
        \fill[matred!30]  (0.0,0.0) rectangle (0.05, 0.05);
        \draw[black, thin] (0.0,0.0) rectangle (0.05, 0.05);
        \node at (0.2,0.02)  {$\mesh^{im}$};
    \end{scope}

    \begin{scope}[shift={(-0.35, 0.05)}]
        \draw[->, thin] (0,0) -- (0.14,0) node[right] {$x$};
        \draw[->, thin] (0,0) -- (0,0.14) node[above] {$y$};
    \end{scope}

    \node at (0.25,-0.1){$\Gamma_{\!i}$};
    \node at (0.75,-0.1){$\Gamma_{\!i}$};  
    \end{scope}

    \begin{scope}[scale=1, xshift=5cm]
\pgfplotsset{every axis/.append style={
                width=8cm,
                height=5cm,
                mark size = 1pt,
                mark options={solid},
                max space between ticks=20,
                semithick,
                xtick pos=bottom,
                xmajorgrids=true, ymajorgrids=true,
                axis line style={thin},
                tick label style = {font=\tiny},
                ytick pos=left,},
        legend style={font=\tiny},
        legend cell align={left},
        legend columns=1,
        transpose legend,
    }

    \begin{groupplot}[group style={group size=1 by 1, horizontal sep=1.5cm, vertical sep=0.5cm}]

        \nextgroupplot[ymode=log, xmode=log, ylabel=$\| \vec{U}_h - \vec{U}_e \|$, xlabel=$\Delta t$, legend to name=mms_legend]

        \addplot[matyellow, mark=square*, dashed, mark options={solid}, forget plot] table [y=L2_U, x=dt, col sep=comma]{L2_errors_imex_1stage_mixed.dat};
        \addplot[matyellow, mark=square*, dashdotted, mark options={solid}, forget plot] table [y=L2_U, x=dt, col sep=comma]{L2_errors_imex_2stage_mixed.dat};
        \addplot[matyellow, mark=square*, mark options={solid}, forget plot] table [y=L2_U, x=dt, col sep=comma]{L2_errors_imex_3stage_mixed.dat};

        \addplot[matorange, mark=square*, dashed, mark options={solid}, forget plot] table [y=L2_U, x=dt, col sep=comma]{L2_errors_imex_1stage_ip.dat};
        \addplot[matorange, mark=square*, dashdotted, mark options={solid}, forget plot] table [y=L2_U, x=dt, col sep=comma]{L2_errors_imex_2stage_ip.dat};
        \addplot[matorange, mark=square*, mark options={solid}, forget plot] table [y=L2_U, x=dt, col sep=comma]{L2_errors_imex_3stage_ip.dat};

        \addplot[matblue, mark=*, dashed, mark options={solid}, forget plot]  table [y=L2_U, x=dt, col sep=comma]{L2_errors_explicit_1stage_ip.dat} [yshift=-5pt] coordinate [pos=0.33] (A) coordinate [pos=0.66] (B);
        \addplot[matblue, mark=*, dashdotted, mark options={solid}, forget plot]  table [y=L2_U, x=dt, col sep=comma]{L2_errors_explicit_2stage_ip.dat} [yshift=12pt]coordinate [pos=0.33] (C) coordinate [pos=0.66] (D);
        \addplot[matblue, mark=*, mark options={solid}, forget plot]  table [y=L2_U, x=dt, col sep=comma]{L2_errors_explicit_3stage_ip.dat} [yshift=-5pt]coordinate [pos=0.33] (E) coordinate [pos=0.66] (F);

        \draw[fill] (B) -| (A) node [xshift=-11pt, yshift=-9pt]  {\tiny{$1$}} node [xshift=2.8pt, yshift=-2.75pt] {\tiny{$1$}} -- cycle;
        \draw[fill] (C) -| (D) node [xshift=11pt, yshift=15pt]  {\tiny{$1$}} node [xshift=-2.8pt, yshift=5.75pt] {\tiny{$2$}} -- cycle;
        \draw[fill] (F) -| (E) node [xshift=-11pt, yshift=-21pt]  {\tiny{$1$}} node [xshift=2.8pt, yshift=-8.75pt] {\tiny{$3$}} -- cycle;

        \addplot[matred, mark=triangle*, dashed, mark options={solid}, forget plot] table [y=L2_U, x=dt, col sep=comma]{L2_errors_implicit_1stage_mixed.dat};
        \addplot[matred, mark=triangle*, dashdotted, mark options={solid}, forget plot] table [y=L2_U, x=dt, col sep=comma]{L2_errors_implicit_2stage_mixed.dat};
        \addplot[matred, mark=triangle*, mark options={solid}, forget plot] table [y=L2_U, x=dt, col sep=comma]{L2_errors_implicit_3stage_mixed.dat};

        \addplot[matyellow, mark=square*, mark options={solid}, draw=none] table [y=L2_U, x=dt, col sep=comma]{L2_errors_imex_3stage_mixed.dat}; \addlegendentry{IMEX-MX};
        \addplot[matorange, mark=square*, mark options={solid}, draw=none] table [y=L2_U, x=dt, col sep=comma]{L2_errors_imex_3stage_ip.dat}; \addlegendentry{IMEX-IP};
        \addplot[matblue, mark=*, mark options={solid},  draw=none] table [y=L2_U, x=dt, col sep=comma]{L2_errors_explicit_3stage_ip.dat}; \addlegendentry{DG-IP};
        \addplot[matred, mark=triangle*, mark options={solid}, draw=none] table [y=L2_U, x=dt, col sep=comma]{L2_errors_implicit_3stage_mixed.dat}; \addlegendentry{HDG-MX};
        
        \addplot[black, dashed, draw=none] table [y=L2_U, x=dt, col sep=comma]{L2_errors_imex_1stage_mixed.dat}; \addlegendentry{RK1};
        \addplot[black, dashdotted, draw=none] table [y=L2_U, x=dt, col sep=comma]{L2_errors_imex_2stage_mixed.dat}; \addlegendentry{RK2};
        \addplot[black, solid, draw=none] table [y=L2_U, x=dt, col sep=comma]{L2_errors_imex_3stage_mixed.dat}; \addlegendentry{RK3};

        \coordinate (top) at (rel axis cs:1,1);
        \coordinate (bottom) at (rel axis cs:1,0);

    \end{groupplot}

    \path (top)--(bottom) coordinate[midway] (group center);
    \node[right=0.5em,inner sep=0pt] at(group center -| current bounding box.east) {\pgfplotslegendfromname{mms_legend}};
    \end{scope}

    \node[above] at (1.5, 3.5) {$(a)$};
    \node[above] at (8.1, 3.5) {$(b)$};

\end{tikzpicture}
    \caption[Method of Manufactured Solutions: Mesh and Errors]{Method of Manufactured Solutions: Representation of the mesh $\mesh$ with element size $h = 0.035$ and interface boundaries $\Gamma_{\!i}$ $(a)$, and corresponding $L^2$ errors of the solution after one period, $t_{\mathrm{end}} = 0.01$, shown over a range of time step sizes $\Delta t$ $(b)$.}
    \label{fig::domain_mms}
\end{figure}

In \Cref{fig::domain_mms} we show the L2 errors of the solution after one period $t_{\textrm{total}} = 2/\omega$ over a range of time step sizes $\Delta t$ for different $s$-stage schemes introduced in \Cref{sec::imex}. As a comparison we included the results of an explicit DG scheme and an implicit HDG-MX scheme on $\mesh$, since they use the same spatial discretization and time integration schemes. The results confirm that the IMEX-MX and IMEX-IP schemes achieve the expected convergence rates.

\subsubsection{Isentropic Vortex in a Uniform Flow}
\label{sec::experiments::vortex_transport_uniform_flow}
\pgfplotsset{table/search path={./figures/vortex_data/}}

This validation benchmark assesses the ability of the proposed IMEX schemes to transport vortical structures in a uniform flow over long time intervals for varying polynomial order~$k$. The domain is $\Omega = [0,1]^2$ with periodic boundary conditions, discretized by a $32 \times 32$ quadrilateral mesh~$\mesh$. The associated implicit mesh $\mesh^{im}$ consists of $8 \times 32$ elements and is refined toward the center, following \cite{vermeire2021110022}.

The initial condition is the fast isentropic vortex from \cite[Section~4.6]{wang2013high}, superimposed on a uniform free-stream velocity $\vec{u}_\infty$. A convective time scale is defined as the time required for the vortex to traverse the domain, $t_c = 1 / |\vec{u}_\infty|$. The simulation is run until $t_{\mathrm{end}} = 5 t_c$, corresponding to five periods. The exact solution corresponds to pure advection of the initial state, $\vec{U}_{e} (\vec{x}, t) = \vec{U}(\vec{x} - \vec{u}_\infty t, 0)$, and errors are measured accordingly. Time integration uses a 3-stage ARS scheme \cite{ascher1997implicit}. Results are compared to a fully explicit DG discretization with a time step restricted by the $k=5$ stability limit, representing the stiffest case.

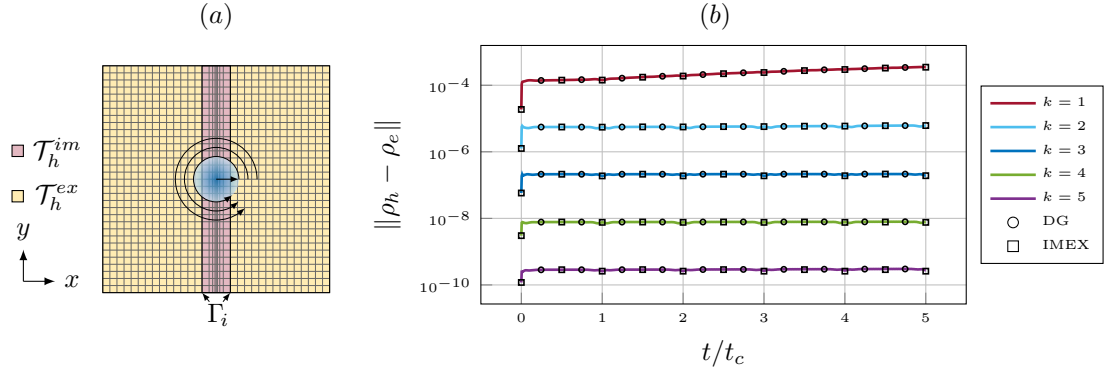
\begin{figure}[ht]
    \centering
    \tikzsetnextfilename{domain_vortex}
\tikzsetnextfilename{domain_vortex}
\begin{tikzpicture}[>=latex]

\def\Xcoords{
    -0.5,-0.46875,-0.4375,-0.40625,-0.375,-0.34375,-0.3125,-0.28125,-0.25,
    -0.21875,-0.1875,-0.15625,-0.125,-0.09375,-0.0625,-0.03125,-0.015625,
    -0.00520833,0,0.00520833,0.015625,0.03125,0.0625,0.09375,0.125,0.15625,
    0.1875,0.21875,0.25,0.28125,0.3125,0.34375,0.375,0.40625,0.4375,0.46875,0.5
}
\def\Ycoords{
    -0.5,-0.46875,-0.4375,-0.40625,-0.375,-0.34375,-0.3125,-0.28125,-0.25,
    -0.21875,-0.1875,-0.15625,-0.125,-0.09375,-0.0625,-0.03125,0,0.03125,
    0.0625,0.09375,0.125,0.15625,0.1875,0.21875,0.25,0.28125,0.3125,0.34375,
    0.375,0.40625,0.4375,0.46875,0.5
}

\begin{scope}[scale=3, yshift=0.05cm]
    
    \fill[amber!30]  (0.0,0.0) rectangle (0.4375,1.0);
    \fill[matred!30]  (0.4375,0.0) rectangle (0.5625,1.0);
    \fill[amber!30]  (0.5625,0.0) rectangle (1.0,1.0);

    \foreach \x in \Xcoords {
        \pgfmathsetmacro{\xx}{\x + 0.5}
        \draw[black!60, ultra thin] (\xx,0) -- (\xx,1);
    }
    \foreach \y in \Ycoords {
        \pgfmathsetmacro{\yy}{\y + 0.5}
        \draw[black!60, ultra thin] (0,\yy) -- (1,\yy);
    }

    \draw[black, thin] (0,0) rectangle (1,1);
    \draw[black, thin] (0.4375,0) rectangle (0.5625,1);

    \begin{scope}[shift={(-0.4, 0.4)}]
        \fill[amber!30]  (0.0,0.0) rectangle (0.05, 0.05);
        \draw[black, thin] (0.0,0.0) rectangle (0.05, 0.05);
        \node at (0.2,0.02)  {$\mesh^{ex}$};
    \end{scope}

    \begin{scope}[shift={(-0.4, 0.6)}]
        \fill[matred!30]  (0.0,0.0) rectangle (0.05, 0.05);
        \draw[black, thin] (0.0,0.0) rectangle (0.05, 0.05);
        \node at (0.2,0.02)  {$\mesh^{im}$};
    \end{scope}

    \begin{scope}[shift={(-0.35, 0.05)}]
        \draw[->, thin] (0,0) -- (0.14,0) node[right] {$x$};
        \draw[->, thin] (0,0) -- (0,0.14) node[above] {$y$};
    \end{scope}

    \node at (0.51,-0.1){$\Gamma_{\!i}$};
    \draw[->, ultra thin] (0.47, -0.04) -- (0.4375, 0);
    \draw[->, ultra thin] (0.53, -0.04) -- (0.5625, 0);

    \shade[outer color=matblue!20, inner color=matblue, opacity=0.7] (0.5, 0.5) circle (0.1);
    
    \draw[-{>[flex=0.85]}, ultra thin] (0.6, 0.5) arc (0:315:0.1);
    \draw[-{>[flex=0.85]}, ultra thin] (0.64, 0.5) arc (0:315:0.14);
    \draw[-{>[flex=0.85]}, ultra thin] (0.68, 0.5) arc (0:315:0.18);
    \draw[-{>[flex=0.85]}, ultra thin] (0.5, 0.5) -- (0.6, 0.5);

\end{scope}

    \begin{scope}[scale=1, xshift=5cm]
\pgfplotsset{every axis/.append style={
                width=8cm,
                height=5cm,
                mark size = 1pt,
                mark repeat = 100,
                semithick,
                xtick pos=bottom,
                xmajorgrids=true, ymajorgrids=true,
                axis line style={thin},
                tick label style = {font=\tiny},
                ytick pos=left,},
        legend style={font=\tiny},
        legend cell align={left},
        legend columns=1,
        transpose legend,
    }

    \begin{groupplot}[group style={group size=1 by 1, horizontal sep=1.5cm, vertical sep=0.5cm}]

        \nextgroupplot[ymode=log, ylabel=$\| \rho_h - \rho_e \|$, xlabel=$t/t_c$, legend to name=vortex_legend, ytick={1e-10, 1e-8, 1e-6, 1e-4}, yticklabels={$10^{-10}$,  $10^{-8}$,  $10^{-6}$, $10^{-4}$}]


        \addplot[matred, line width = 1pt]  table [y=rho, x=t, col sep=comma]{dg_1_rho.csv}; \addlegendentry{$k=1$};
        \addplot[matazure, line width = 1pt]  table [y=rho, x=t, col sep=comma]{dg_2_rho.csv};\addlegendentry{$k=2$};
        \addplot[matblue, line width = 1pt]  table [y=rho, x=t, col sep=comma]{dg_3_rho.csv};\addlegendentry{$k=3$};
        \addplot[matgreen, line width = 1pt]  table [y=rho, x=t, col sep=comma]{dg_4_rho.csv};\addlegendentry{$k=4$};
        \addplot[matpurple, line width = 1pt]  table [y=rho, x=t, col sep=comma]{dg_5_rho.csv};\addlegendentry{$k=5$};

        \addplot[mark=o, mark phase=50, only marks]  table [y=rho, x=t, col sep=comma]{dg_1_rho.csv}; \addlegendentry{DG};
        \addplot[mark=o, mark phase=50, only marks, forget plot]  table [y=rho, x=t, col sep=comma]{dg_2_rho.csv};
        \addplot[mark=o, mark phase=50, only marks, forget plot]  table [y=rho, x=t, col sep=comma]{dg_3_rho.csv};
        \addplot[mark=o, mark phase=50, only marks, forget plot]  table [y=rho, x=t, col sep=comma]{dg_4_rho.csv};
        \addplot[mark=o, mark phase=50, only marks, forget plot]  table [y=rho, x=t, col sep=comma]{dg_5_rho.csv};

        \addplot[mark = square, mark phase=0, only marks] table [y=rho, x=t, col sep=comma]{imex_1_rho.csv}; \addlegendentry{IMEX};
        \addplot[mark = square, mark phase=0, only marks] table [y=rho, x=t, col sep=comma]{imex_2_rho.csv};
        \addplot[mark = square, mark phase=0, only marks] table [y=rho, x=t, col sep=comma]{imex_3_rho.csv};
        \addplot[mark = square, mark phase=0, only marks] table [y=rho, x=t, col sep=comma]{imex_4_rho.csv};
        \addplot[mark = square, mark phase=0, only marks] table [y=rho, x=t, col sep=comma]{imex_5_rho.csv};


        \coordinate (top) at (rel axis cs:1,1);
        \coordinate (bottom) at (rel axis cs:1,0);

    \end{groupplot}

    \path (top)--(bottom) coordinate[midway] (group center);
    \node[right=0.5em,inner sep=0pt] at(group center -| current bounding box.east) {\pgfplotslegendfromname{vortex_legend}};
    \end{scope}

    \node[above] at (1.5, 3.5) {$(a)$};
    \node[above] at (8.1, 3.5) {$(b)$};

\end{tikzpicture}
    \caption[Isentropic Vortex Transport by Uniform Flow: Mesh and Errors]{Isentropic Vortex Transport by Uniform Flow: Representation of the mesh $\mesh$ with refined elements and interface boundaries $\Gamma_{\!i}$ $(a)$, and corresponding $L^2$ errors of the solution after $5$ periods $(b)$.}
    \label{fig::domain_vortex}
\end{figure}

\Cref{fig::domain_vortex} shows the $L^2$ density error versus $t/t_c$. The IMEX scheme\footnote{In the absence of viscous terms, IMEX-MX and IMEX-IP coincide.} maintains bounded errors throughout the simulation\footnote{For $k=1$, the apparent growth is caused by excessive numerical diffusion, which gradually smears the vortex.}. Accuracy is comparable to the explicit DG scheme, indicating that the IMEX approach preserves fidelity in convection-dominated regimes. No error increase is observed as the vortex crosses the interface boundaries~$\Gamma_{\!i}$, confirming the consistency of the interface treatment. The observed errors agree well with theoretical expectations for high-order discretizations of smooth solutions.

\subsubsection{Laminar Flow over a NACA-0012 Airfoil}
\label{sec::experiments::naca0012_airfoil}
\pgfplotsset{table/search path={./figures/naca_data/}}

\begin{figure}
    \centering
    \includegraphics{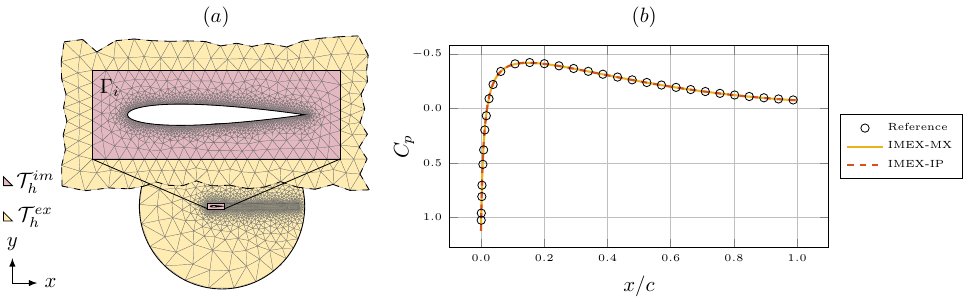}
    \caption[Laminar Flow over a NACA-0012 Airfoil: Mesh and Pressure Coefficient]{Laminar Flow over a NACA-0012 Airfoil: Implicit (red) and explicit (yellow) meshes with interface boundary $\Gamma_{\!i}$ over the NACA-0012 $(a)$ and 
    the pressure coefficient $C_p$ over the airfoil surface; The reference $C_p$ values are digitized \cite{WebPlotDigitizer} from [Figure~13a, \citenum{swanson2016comparison}].}
    \label{fig::naca0012_domain}
\end{figure}

Our second validation case is based on the full Navier-Stokes equations. We consider the classical problem of two-dimensional laminar flow around a NACA-0012 airfoil \cite{caraeni2002, venkatakrishnan1990191, swanson2016comparison} at zero angle of attack. The flow is characterized by a Reynolds number of $\mathrm{\Re}_c = 5000$, defined with respect to the chord length $c = 1$ (meters) and free-stream conditions ($\rho_\infty$, $u_\infty$, $\mu_\infty$).

The computational domain extends radially to a distance of $7c$ from the trailing edge. On the airfoil surface, adiabatic wall boundary conditions are imposed. The outer boundary is treated using far-field conditions in combination with a sponge layer of width $1.5c$ to minimize reflections.

An unstructured triangular mesh is employed throughout both the implicit and explicit regions. Although a structured quadrilateral mesh (e.g., an O-grid) is generally preferable for accurately resolving boundary layers, solver limitations prevent the generation of a fully body-fitted structured grid. To mitigate this, an O-grid topology (rounded near the trailing edge) is used within the implicit region to better align elements with the boundary layer, see \Cref{fig::naca0012_domain}. The near-wall resolution corresponds to a wall-normal spacing of $\Delta y \approx 0.00613c$, with comparable streamwise spacing ($\Delta x \approx \Delta y$) maintained around the airfoil\footnote{These spacings refer to the characteristic element size. Owing to the use of high-order polynomial approximations, the effective resolution based on the solution degrees of freedom is higher.}.

The flow is initialized with free-stream conditions and advanced in time up to $t_f = 30 t_c$, where the convective time scale is defined as $t_c = c / u_\infty$, representing the time required for the flow to traverse one chord length. \Cref{fig::naca0012_contours} shows an instantaneous vorticity $\omega$ and pressure $p$ contour for the IMEX-MX schemes.

The computed lift and drag coefficients show good agreement with reference data. Specifically, the lift coefficients converge to $C_l^{mx} \approx 6.12 \times 10^{-5}$ and $C_l^{ip} \approx 1.81 \times 10^{-4}$, for IMEX-MX and IMEX-IP, respectively, while the drag coefficients are $C_d^{mx} \approx 0.055$ and $C_d^{ip} \approx 0.0549$. These results are consistent with the reference values reported in \cite{swanson2016comparison}, where $C_l \approx 0$ and $C_d \approx 0.0556$. Furthermore, the pressure coefficient $C_p$ in \Cref{fig::naca0012_domain} (right panel) for both IMEX methods shows excellent agreement with the same reference data. The estimated separation point, based on the wall shear stress, is approximately $x_{sep} \approx 0.807c$ for both IMEX-IP and IMEX-MX, also in close agreement with the literature.

\subsection{Computational Performance Assessment}

Having established the correctness of the proposed solvers through a series of verification and validation benchmarks, we now turn our attention to evaluating their computational performance. To ensure consistency, all experiments are conducted on a fixed hardware configuration, consisting of a single compute node equipped with two \texttt{AMD EPYC 7713} CPUs and $512$ GB of RAM.

The present study is restricted to two-dimensional benchmark configurations, following, e.g., \cite{pereira2025113819,franciolini2020104542}. This setting provides a controlled and well-understood environment in which the performance characteristics of the proposed methods can be assessed in isolation. The resulting observations therefore serve as a baseline for the guidelines developed in this work. Although these configurations correspond to laminar flow regimes, the underlying trends are expected to become even more pronounced in three-dimensional simulations, where additional complexity arising from turbulence (particularly in the context of WR-LES) further amplifies stiffness and computational cost.

The linear systems arising within the implicit nonlinear solver are solved using the \texttt{PARDISO} sparse direct solver from the \texttt{Intel MKL} library \cite{schenk2004475,intel_mkl}. Throughout this study, convergence is declared when the residual is reduced below $10^{-14}$. 
 
\begin{remark}
	While sparse direct solvers are not well-suited for large-scale problems due to their computational and memory costs, their use here enables precise measurement of the implicit wall-clock time, $\mathcal{C}^{\mathrm{im}}$, and provides a reliable baseline for future comparisons with more scalable iterative solvers.
\end{remark}

\subsubsection{Speedup - Laminar Flow over a Cylinder}
\label{subsec::cylinder_benchmark}
\pgfplotsset{table/search path={./figures/cylinder_data/}}

This benchmark considers a two-dimensional flow past a circular cylinder in a configuration similar to that of \cite{pereira2025113819}. Instead of assessing the speedup $S$ of a single spatial resolution, we examine a range of meshes and varying numbers of implicitly treated elements. This approach provides a more comprehensive characterization of the performance of the IMEX scheme across different resolutions and implicit-explicit meshes. It also enables identification of general trends in how IMEX efficiency varies with boundary-layer meshing strategies, providing practical guidance on expected performance.

The computational domain is defined as the annular region
\begin{align}
	\Omega = \left\{ (x,y) \;|\; R \le \sqrt{x^2+y^2} \le 100 R \right\},
\end{align}
where $R = 0.5$ denotes the cylinder radius. Adiabatic no-slip wall boundary conditions are imposed on the cylinder surface $\Gamma_{\!w}$, while far-field boundary conditions are prescribed on the outer boundary $\Gamma_\infty$. Note that units are omitted for readability; all dimensional quantities are expressed in SI units.

The flow regime is laminar, characterized by $\Ma_\infty = 0.2$ and $\Re_\infty = 150$, based on the free-stream conditions and cylinder diameter. In the wall-normal (radial) direction, an estimated first-element size of $h_0 = 0.05$ is adopted \cite{inoueSoundGenerationTwodimensional2002}, consistent with the laminar boundary layer scaling $\mathcal{O}( \Re^{-1/2})$. To systematically investigate the influence of geometric stiffness in the boundary layer, we consider a set of meshes of the form
\begin{alignat*}{3}
	\mesh(n_{i}, h_{i}) &:= \mesh^{im}(n_{i}, h_{i}) \cup \mesh^{ex}(n - n_{i}, h_0), & \qquad
	h_i &:= \min\limits_{T \in \mesh^{im}}{h_T}.
\end{alignat*}

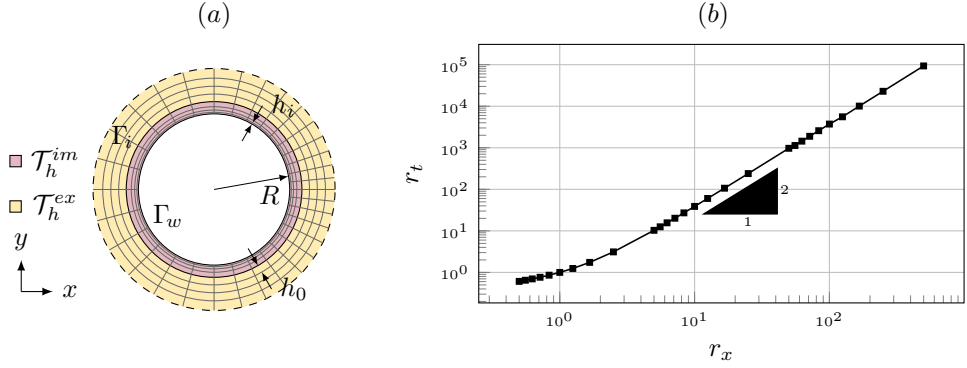
\begin{figure}
    \centering
    \tikzsetnextfilename{domain_cylinder}
\tikzsetnextfilename{domain_cylinder}
\begin{tikzpicture}[>=latex]






\begin{scope}[scale=2, shift={(0.75, 0.75)}]

    \fill[matred!30, even odd rule] (0.0, 0.0) circle (0.5) (0.0, 0.0) circle (0.58075118);
    \fill[amber!30, even odd rule] (0.0, 0.0) circle (0.58075118) (0.0, 0.0) circle (0.8);
\draw[black, thin] (0.0, 0.0) circle (0.5);
\draw[black, thin] (0.0, 0.0) circle (0.58075118);
\draw[black, thin, dashed] (0.0, 0.0) circle (0.8);

\foreach \r in {0.51, 0.52495349,  0.54731417, 0.63075118, 0.68249636, 0.73604762} {
  \draw [black!60, ultra thin] (0.0, 0.0) circle (\r);
}

\foreach \phi in { 0.        ,   5.625     ,  11.81709606,  18.63346113,
        26.13703221,  34.39709139,  43.48990559,  53.49943073,
        64.51808693,  76.64761185,  90.        , 105.        ,
       120.        , 135.        , 150.        , 165.        ,
       180.        , 195.        , 210.        , 225.        ,
       240.        , 255.        , 270.        , 283.35238815,
       295.48191307, 306.50056927, 316.51009441, 325.60290861,
       333.86296779, 341.36653887, 348.18290394, 354.375} {

  \draw [black!60, ultra thin, shift={(0.0, 0.0)}] ({0.5 * cos(\phi)}, {0.5 * sin(\phi)}) -- ({0.8 * cos(\phi)}, {0.8 * sin(\phi)});
}

\def\angh0{60}
\draw[very thin, ->] ({0.37 * cos(\angh0)}, {0.37 * sin(\angh0)}) -- ({0.5 * cos(\angh0)}, {0.5 * sin(\angh0)});
\draw[very thin, ->] ({0.64 * cos(\angh0)}, {0.64 * sin(\angh0)}) -- ({0.51 * cos(\angh0)}, {0.51 * sin(\angh0)});
\draw[very thin, draw=none] ({0.37 * cos(\angh0)}, {0.37 * sin(\angh0)}) -- ({0.64 * cos(\angh0)}, {0.64 * sin(\angh0)}) node[anchor=west] {$h_i$};

\def\angh1{300}
\draw[very thin, ->] ({0.45 * cos(\angh1)}, {0.45 * sin(\angh1)}) -- ({0.58075118 * cos(\angh1)}, {0.58075118 * sin(\angh1)});
\draw[very thin, ->] ({0.76 * cos(\angh1)}, {0.76 * sin(\angh1)}) -- ({0.63075118 * cos(\angh1)}, {0.63075118 * sin(\angh1)});
\draw[very thin, draw=none] ({0.45 * cos(\angh1)}, {0.45 * sin(\angh1)}) -- ({0.76 * cos(\angh1)}, {0.76 * sin(\angh1)}) node[anchor=west] {$h_0$};

\def\angh4{10}
\draw[very thin, ->] (0.0, 0.0) -- ({0.5 * cos(\angh4)}, {0.5 * sin(\angh4)}) node[anchor=north east] {$R$};

\def\angh2{210}
\node at ({0.35 * cos(\angh2)}, {0.35 * sin(\angh2)}) {$\Gamma_{\!w}$};

\def\angh3{150}
\node at ({0.7 * cos(\angh3)}, {0.7 * sin(\angh3)}) {$\Gamma_{\!i}$};
\end{scope}

\begin{scope}[scale=3]
\begin{scope}[shift={(-0.4, 0.4)}]
    \fill[amber!30]  (0.0,0.0) rectangle (0.05, 0.05);
    \draw[black, thin] (0.0,0.0) rectangle (0.05, 0.05);
    \node at (0.2,0.02)  {$\mesh^{ex}$};
\end{scope}

\begin{scope}[shift={(-0.4, 0.6)}]
    \fill[matred!30]  (0.0,0.0) rectangle (0.05, 0.05);
    \draw[black, thin] (0.0,0.0) rectangle (0.05, 0.05);
    \node at (0.2,0.02)  {$\mesh^{im}$};
\end{scope}

\begin{scope}[shift={(-0.35, 0.05)}]
    \draw[->, thin] (0,0) -- (0.14,0) node[right] {$x$};
    \draw[->, thin] (0,0) -- (0,0.14) node[above] {$y$};
\end{scope}
\end{scope}



\begin{scope}[scale=1, xshift=5cm]
\pgfplotsset{every axis/.append style={
                width=8cm,
                height=5cm,
                mark size = 1pt,
                mark options={solid},
                max space between ticks=20,
                semithick,
                xtick pos=bottom,
                xmajorgrids=true, ymajorgrids=true,
                axis line style={thin},
                tick label style = {font=\tiny},
                ytick pos=left,},
        legend style={font=\tiny},
        legend cell align={left},
        legend columns=1,
        transpose legend,
    }

    \begin{groupplot}[group style={group size=1 by 1, horizontal sep=1.5cm, vertical sep=0.5cm}]

        \nextgroupplot[ymode=log, xmode=log, ylabel=$r_t$, xlabel=$r_x$]

        \addplot[black, mark=square*] table [y=r_t, x=r_x, col sep=comma]{stable_explicit_dt.dat} [yshift=-5pt] coordinate [pos=0.4] (A) coordinate [pos=0.6] (B);
        \draw[fill] (B) -| (A) node [xshift=-11pt, yshift=-20pt]  {\tiny{$1$}} node [xshift=2.8pt, yshift=-8pt] {\tiny{$2$}} -- cycle;

    \end{groupplot}
    \end{scope}

    \node[above] at (1.5, 3.5) {$(a)$};
    \node[above] at (8.1, 3.5) {$(b)$};

\end{tikzpicture}
    \caption[Speedup - Laminar Flow over a Cylinder: Mesh and Stability]{Speedup - Laminar Flow over a Cylinder: Representative mesh $\mesh(n_{i} = 128, h_i = 0.01)$ showing interface boundaries $\Gamma_{\!i}$ $(a)$, and the numerically computed stable time-step ratio $r_t$ as a function of geometric stiffness $r_x$ $(b)$.}
    \label{fig::2_domain}
\end{figure}

The total number of quadrilateral elements is fixed at $n = 4096$, corresponding to a $128 \times 32$ discretization in the radial and azimuthal directions, respectively, with the azimuthal resolution always held constant at $32$ elements. 
By \textit{fixing the smallest element size in the explicit mesh} to $h_0$ and progressively refining the implicit region such that $h_i < h_0$, we observe the following:
\begin{itemize}
	\item The geometric stiffness, defined as $r_x = h_0 / h_i$, increases. As a consequence, the stability constraint of the fully explicit scheme becomes more restrictive, leading to a reduction in the stable time step~$\Delta t_s^{ex}(\mesh, \vec{U}_h)$.
	\item In contrast, the stable time step of the IMEX schemes, $\Delta t_s^{imex}(\mesh^{ex}, \vec{U}_h)$, remains unchanged, as it is governed solely by the fixed explicit element size $h_0$. Consequently, the ratio of stable time steps $r_t$, increases with increasing geometric stiffness.
	\item The number of nonlinear iterations required during the implicit stage remains constant at $m = 2$, see \Cref{eq::implicit_wall_clock_time}, indicating that the cost of the implicit solve is unaffected by the refinement of $\mesh^{im}$.
\end{itemize}

This setup enables a systematic assessment of the speedup of the IMEX scheme relative to its fully explicit DG counterpart, as a function of both the number of implicit elements $n_{i}$ and the geometric stiffness $r_x$. In \Cref{fig::2_domain} we show a representative mesh $\mesh$ with $n_{i} = 128$
implicit elements and an element size of $h_i = 0.01$, along with the relationship between the stable time-step ratio $r_t$ and the geometric stiffness $r_x$. The values of $r_t$ are obtained via a bisection procedure~\cite{pereira2025113819}, in which the explicit time step $\Delta t^{ex}$ is iteratively adjusted until the stability threshold is identified. 

The flow is initialized with a uniform free-stream state and advanced in time until transient effects decay and a periodic vortex-shedding regime is established. The resulting Strouhal number $\St \approx 0.18$, the lift and drag coefficients (see \Cref{fig::cylinder_coeffcients}), serve as a reference indicator of physical consistency \cite{inoueSoundGenerationTwodimensional2002}. Once this regime is reached, the wall-clock time per time step is measured 
over one shedding period $t_{\text{total}} \approx 1/\St$ for simulations run via an IMEX and a fully explicit schemes on the same mesh $\mesh(n_{i}, h_i)$. \Cref{fig::cylinder_domain_runtimes} indicates that the measurements are not dominated by noise, 
hence the mean wall-clock times are obtained by averaging over the shedding period and used to compute the speedup $S$ as defined in \Cref{eq::speedup}.

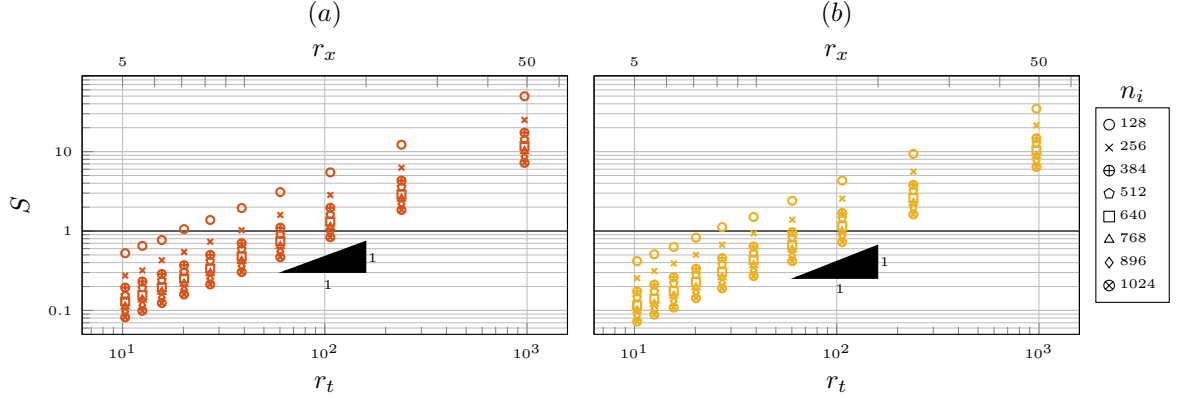
\begin{figure}
    \centering
    \tikzsetnextfilename{linear_scaling}
\tikzsetnextfilename{linear_scaling}
\begin{tikzpicture}
    \pgfplotsset{every axis/.append style={
        width=8cm, 
        height=5cm, 
        mark size = 1.5pt,
        mark options={solid},
         max space between ticks=20,
         thick,
         xtick pos=bottom,
         xmajorgrids=true, ymajorgrids=true,
         yminorgrids=true,
         axis line style={thin}, 
         tick label style = {font=\tiny},
         ytick pos=left,
        },        
        legend style={font=\tiny},
        legend cell align={left},
        legend columns=1,
        }

        \begin{groupplot}[group style={group size=2 by 2,  vertical sep=0.5cm, horizontal sep=10pt}]
            
            \nextgroupplot[xmode=log, ymode=log, ylabel = $S$, xlabel = $r_t$, yticklabels={$0.01$, $0.1$, $1$, $10$, $100$}, extra y ticks={1}, extra y tick labels = {}, extra y tick style={ grid=major, grid style={black}, tick style={draw=none}}, ymin=0.05, ymax=90, legend to name=legend_speedup_linear,]
            
            \addlegendimage{only marks, mark=o};        \addlegendentry{$128$};
            \addlegendimage{only marks, mark=x};        \addlegendentry{$256$};
            \addlegendimage{only marks, mark=oplus};    \addlegendentry{$384$};
            \addlegendimage{only marks, mark=pentagon};   \addlegendentry{$512$};
            \addlegendimage{only marks, mark=square};   \addlegendentry{$640$};
            \addlegendimage{only marks, mark=triangle}; \addlegendentry{$768$};
            \addlegendimage{only marks, mark=diamond};  \addlegendentry{$896$};
            \addlegendimage{only marks, mark=otimes}; \addlegendentry{$1024$};

            \addplot[only marks, matorange, mark=o]         table [y=4, x=tau, col sep=comma]{speedup_ip.dat};      
            \addplot[only marks, matorange, mark=x]         table [y=8, x=tau, col sep=comma]{speedup_ip.dat};      
            \addplot[only marks, matorange, mark=oplus]     table [y=12, x=tau, col sep=comma]{speedup_ip.dat};      
            \addplot[only marks, matorange, mark=pentagon]    table [y=16, x=tau, col sep=comma]{speedup_ip.dat};      
            \addplot[only marks, matorange, mark=square]    table [y=20, x=tau, col sep=comma]{speedup_ip.dat};      
            \addplot[only marks, matorange, mark=triangle]  table [y=24, x=tau, col sep=comma]{speedup_ip.dat};      
            \addplot[only marks, matorange, mark=diamond]  table [y=28, x=tau, col sep=comma]{speedup_ip.dat};      
            \addplot[only marks, matorange, mark=otimes]  table [y=32, x=tau, col sep=comma]{speedup_ip.dat};      
            \addplot[black, dashed, thick, domain=10:1000, samples=2,forget plot, draw=none] { 5.92/990 * (x-10) + 0.08}[yshift=-5pt] coordinate [pos=0.4] (A) coordinate [pos=0.6] (B);
            \draw[fill] (A) -| (B) node [xshift=-14pt, yshift=-16pt]  {\tiny{$1$}} node [xshift=2.8pt, yshift=-6pt] {\tiny{$1$}} -- cycle;

            \path (rel axis cs:0,1)--(rel axis cs:1,1) coordinate[midway] (group center_left);

            \nextgroupplot[xmode=log, ymode=log, xlabel = $r_t$, y tick style={draw=none}, yticklabels={}, extra y tick labels = {}, extra y ticks={1}, extra y tick style={ grid=major, grid style={black}, tick style={draw=none}}, ymin=0.05, ymax=90]
            
            \addplot[only marks, matyellow, mark=o]        table [y=4, x=tau, col sep=comma]{speedup_mixed.dat};      
            \addplot[only marks, matyellow, mark=x]        table [y=8, x=tau, col sep=comma]{speedup_mixed.dat};      
            \addplot[only marks, matyellow, mark=oplus]    table [y=12, x=tau, col sep=comma]{speedup_mixed.dat};      
            \addplot[only marks, matyellow, mark=pentagon]     table [y=16, x=tau, col sep=comma]{speedup_mixed.dat};      
            \addplot[only marks, matyellow, mark=square]   table [y=20, x=tau, col sep=comma]{speedup_mixed.dat};      
            \addplot[only marks, matyellow, mark=triangle]   table [y=24, x=tau, col sep=comma]{speedup_mixed.dat};      
            \addplot[only marks, matyellow, mark=diamond]    table [y=28, x=tau, col sep=comma]{speedup_mixed.dat};      
            \addplot[only marks, matyellow, mark=otimes] table [y=32, x=tau, col sep=comma]{speedup_mixed.dat};      
            \addplot[black, dashed, thick, domain=10:1000, samples=2,forget plot, draw=none] { 5.94/990 * (x-10) + 0.06} [yshift=-5pt] coordinate [pos=0.4] (A) coordinate [pos=0.6] (B);
            \draw[fill] (A) -| (B) node [xshift=-14pt, yshift=-16pt]  {\tiny{$1$}} node [xshift=2.8pt, yshift=-6pt] {\tiny{$1$}} -- cycle;

            \coordinate (right_bottom) at (rel axis cs:1,0);
            \coordinate (right_top) at (rel axis cs:1,1);
            \path (rel axis cs:0,1)--(rel axis cs:1,1) coordinate[midway] (group center_right);
            
            \nextgroupplot[axis x line*=top, xmode=log, ymode=log, at=(group c1r1.north), xtick={5,6,7,8,9,10,20,30,40,50, 60}, xticklabels={},
            extra x tick labels = {$5$, $50$}, extra x ticks={5, 50}, extra x tick style={tick style={draw=none}},
            y tick style={draw=none}, yticklabels={},  xlabel=$r_x$, xlabel shift=-8pt, grid=none]
            \addplot[black, dashed, thick, domain=5:50, samples=2, draw=none] { 0.54/45* (x-5) + 0.06};

            \nextgroupplot[axis x line*=top, xmode=log, ymode=log, at=(group c2r1.north), xtick={5,6,7,8,9,10,20,30,40,50, 60}, xticklabels={},
            extra x tick labels = {$5$, $50$}, extra x ticks={5, 50}, extra x tick style={tick style={draw=none}},
            y tick style={draw=none}, yticklabels={},  xlabel=$r_x$, xlabel shift=-8pt, grid=none]
            \addplot[black, dashed, thick, domain=5:50, samples=2, draw=none] { 0.54/45* (x-5) + 0.06};

        \end{groupplot}
        \path (right_bottom)--(right_top) coordinate[midway] (group center);
        \node[inner sep=0pt,xshift=20pt] at(group center) {\pgfplotslegendfromname{legend_speedup_linear}};
        \node[inner sep=0pt,xshift=20pt, yshift=-7pt] at(right_top) {$n_i$};
        \node[inner sep=0pt, yshift=23pt] at(group center_left) {$(a)$};
        \node[inner sep=0pt, yshift=23pt] at(group center_right) {$(b)$};

\end{tikzpicture}
    \caption[Speedup - Laminar Flow over a Cylinder: Speedup for IMEX-IP and IMEX-MX]{Speedup - Laminar Flow over a Cylinder: Speedup for different meshes $\mesh(n_{i}, h_i)$ for IMEX-IP $(a)$ and IMEX-MX $(b)$ as a function the stable time step ratio $r_t$ (bottom axis) and geometric stiffness $r_x$ (top axis).}
    \label{fig::linear_scaling_results}
\end{figure}

The results for both IMEX-IP and IMEX-MX are presented in \Cref{fig::linear_scaling_results}. Not surprisingly, the computational cost of the IMEX schemes is dominated by the implicit region, where the nonlinear solves are performed. Consequently, the overall efficiency is strongly influenced by both the size of the implicit mesh and the associated stability constraints. A clear increase in speedup is observed with increasing $r_t$, with speedups reaching up to $S \approx 50$ for IP-based IMEX schemes at $r_x \approx 50$. Moreover, IMEX-IP formulations (HDG-IP and DG-IP) consistently outperform their IMEX-MX counterparts (HDG-MX and DG-IP), a trend that will be examined further in the next benchmark.

However, this behavior does not persist across all configurations. As $r_x$ and $r_t$ decrease (particularly as the number of implicit elements $n_{i}$ increases), the speedup $S$ deteriorates and may fall below that of the fully explicit scheme. This reflects the increasing dominance of the implicit solve cost, which can outweigh the benefit of larger stable time steps. Therefore, \textit{careful partitioning of $\mesh^{im}$ and $\mesh^{ex}$ is essential; otherwise, the IMEX strategy may even lead to a net slowdown.}

\subsubsection[Implicit wall clock time]{Computational Cost of the Implicit Solver}
\pgfplotsset{table/search path={./figures/pulse_data/}}

Having demonstrated the expected IMEX speedups $S$ on a curved mesh in \Cref{subsec::cylinder_benchmark}, and observed the superior performance of IMEX-IP over IMEX-MX, we now examine the underlying causes of this behavior. To this end, we consider a simplified, non-curved configuration in order to isolate the performance characteristics of the implicit solver. We decompose its cost into the principal components, namely system assembly and linear system solution, and analyze how these contributions scale with respect to the number of elements, polynomial degree and the level of shared-memory parallelism (i.e., number of threads). We also illustrate the benefits of having HDG methods over standard DG in implicit systems.

In our implementation, each implicit evaluation comprises the following steps:
\begin{enumerate}
	\item \textbf{Assembly of the linear system}, including the computation of the Jacobian matrix and residual vector. For HDG schemes, this additionally involves forming the Schur complement in \Cref{eq::SchurCompliment}, which requires the dense \textit{local} inversion of $\mathbb{K}$.
	\item \textbf{Solution of the linear system}, consisting of the factorization and back-substitution of the Jacobian matrix, along with the solution of the local problems defined in \Cref{eq::local_solve_static_condensation}.
\end{enumerate}

The mesh $\mesh$ consists of $n_{i}$ quadrilateral elements with periodic boundary conditions on the computational domain $\Omega = [0, \sqrt{n_{i}}]^2$. As an initial condition, we consider a linear (quasi-)one-dimensional acoustic pulse superimposed on a uniform free-stream state, following [Section~5.1, \citenum{ellmenreichCharacteristicBoundaryConditions2026}]. The only modification is the inclusion of viscous effects, allowing for a consistent comparison between the HDG-IP and HDG-MX schemes. The mean wall-clock times for assembly and solution, denoted by $\overline{\mathcal{C}}_{\mathrm{assem}}(n_{i})$ and $\overline{\mathcal{C}}_{\mathrm{solve}}(n_{i})$, are obtained by averaging over $m=500$ iterations in order to reduce statistical noise, see \Cref{fig::pulse_domain_runtimes}.
\begin{figure}[hb]
    \centering
    \tikzsetnextfilename{implicit_timings_over_p}
\tikzsetnextfilename{implicit_timings_per_p}

\begin{tikzpicture}
\pgfplotsset{every axis/.append style={
    ymajorgrids,
    bar width=4pt,
    width=8cm,
    xtick=data,
    width=8cm, 
    height=5cm, 
    tick label style = {font=\tiny},
    legend style={font=\tiny},
    legend cell align={left},
    legend columns=1,
}
}
\begin{groupplot}[group style={group size=2 by 1, x descriptions at=edge bottom, horizontal sep=35pt}, ]

    \nextgroupplot[xlabel=${l}$, title=${(a)}$ ,  ylabel=$\overline{\mathcal{C}}_{\mathrm{assem}}^{\mathrm{MX}}\, /\, \overline{\mathcal{C}}_{\mathrm{assem}}^{\mathrm{IP}}$, symbolic x coords={4, 8, 16, 24, 32, 48, 64}, xticklabels={$4$, $8$, $16$, $24$, $32$, $48$, $64$}, ymin=-0.5, ymax=3.5,legend to name=legend_bar_plot_k, legend image code/.code={
        \draw [#1] (0cm,-0.1cm) rectangle (0.2cm,0.25cm); },]

    \addlegendimage{fill=matred};        
    \addlegendimage{fill=matazure};   
    \addlegendimage{fill=matblue};        
    \addlegendimage{fill=matgreen};        
    \addlegendimage{fill=matpurple};        
    \legend{$1$, $2$, $3$, $4$, $5$}

    \addplot[ybar, bar shift = -8pt, fill=matred] table [y=assemble,x=threads, col sep=comma] {ratio_1_threads.csv};
    \addplot[ybar, bar shift = -4pt, fill=matazure] table [y=assemble,x=threads, col sep=comma] {ratio_2_threads.csv};
    \addplot[ybar, bar shift = -0pt, fill=matblue] table [y=assemble,x=threads, col sep=comma] {ratio_3_threads.csv}; 
    \addplot[ybar, bar shift = 4pt , fill=matgreen] table [y=assemble,x=threads, col sep=comma]  {ratio_4_threads.csv};
    \addplot[ybar, bar shift = 8pt , fill=matpurple] table [y=assemble,x=threads, col sep=comma]  {ratio_5_threads.csv};

    \nextgroupplot[xlabel=${l}$, title=${(b)}$, ylabel=$\overline{\mathcal{C}}_{\mathrm{solve}}^{\mathrm{MX}}\, /\, \overline{\mathcal{C}}_{\mathrm{solve}}^{\mathrm{IP}}$, symbolic x coords={4, 8, 16, 24, 32, 48, 64}, xticklabels={$4$, $8$, $16$, $24$, $32$, $48$, $64$}, ymin=-0.5, ymax=3.5, yticklabels = {}]

    \addplot[ybar, bar shift = -8pt, fill=matred] table [y=solve,x=threads, col sep=comma]  {ratio_1_threads.csv};
    \addplot[ybar, bar shift = -4pt, fill=matazure] table [y=solve,x=threads, col sep=comma]  {ratio_2_threads.csv};
    \addplot[ybar, bar shift = -0pt, fill=matblue] table [y=solve,x=threads, col sep=comma]  {ratio_3_threads.csv}; 
    \addplot[ybar, bar shift = 4pt, fill=matgreen] table  [y=solve,x=threads, col sep=comma]  {ratio_4_threads.csv};
    \addplot[ybar, bar shift = 8pt, fill=matpurple] table  [y=solve,x=threads, col sep=comma]  {ratio_5_threads.csv};

    \coordinate (right_bottom) at (rel axis cs:1,0);
    \coordinate (right_top) at (rel axis cs:1,1);

\end{groupplot}

\path (right_bottom)--(right_top) coordinate[midway] (group center);
\node[inner sep=0pt,xshift=20pt] at(group center) {\pgfplotslegendfromname{legend_bar_plot_k}};
\node[inner sep=0pt,xshift=20pt, yshift=-7pt] at(right_top) {$k$};
\end{tikzpicture}
    \caption[Computational Cost of the Implicit Solver: Assembly and Solve 1]{Computational Cost of the Implicit Solver: Averaged wall-clock times for the acoustic pulse benchmark using HDG-MX (normalized by HDG-IP) for assembly $(a)$ and solve $(b)$, as functions of the number of threads $l$ and polynomial order $k$, with a total of $n_{i} = 128^2$ elements.}
    \label{fig::timings_over_p}
\end{figure}
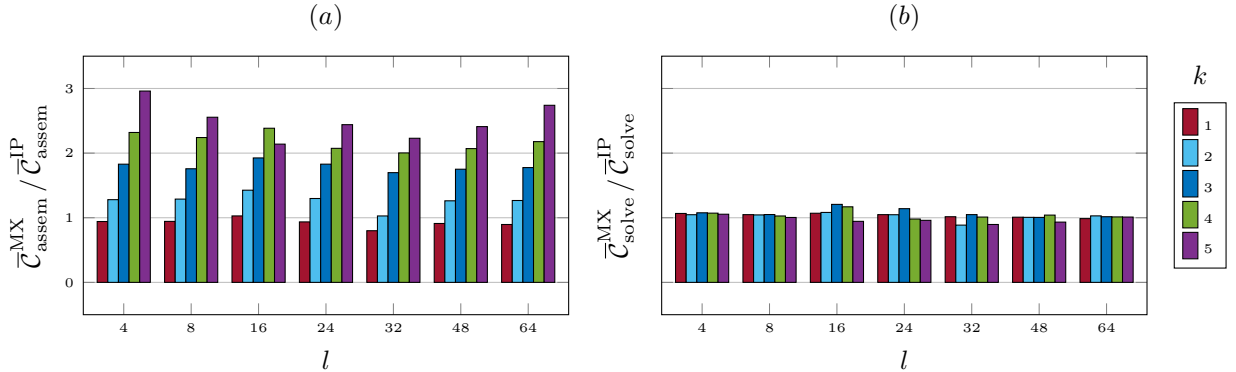


In \Cref{fig::timings_over_p}, we report the averaged wall-clock time (per iteration) for implicit assembly, $\overline{\mathcal{C}}_{\mathrm{assem}}(n_{i})$, and solution, $\overline{\mathcal{C}}_{\mathrm{solve}}(n_{i})$, for an HDG-MX (normalized by HDG-IP) on a mesh with $n_{i} = 128^2$ elements. The results show that the assembly cost of HDG-MX increases markedly with polynomial order $k$, compared to HDG-IP, reaching a slowdown of more than a factor of two for $k = 5$. This increase is attributed to the additional auxiliary terms in the Jacobian matrix that must be assembled and locally inverted in the HDG-MX formulation, see \Cref{rem::hdg_cost}. As a matter of fact, these findings motivate the use of a second layer of static condensation to reduce the size of the local systems \cite{kronbichlerComparisonImplicitExplicit2016}, particularly for high-order discretizations and three-dimensional simulations.

In contrast, the solution times, $\overline{\mathcal{C}}_{\mathrm{solve}}$ of both methods remain comparable across all polynomial orders. This is expected, as both formulations yield linear systems of \textit{identical} sparsity structure, arising from the Schur complement matrices in \Cref{eq::SchurCompliment}.

\begin{figure}[hb]
    \centering
    \tikzsetnextfilename{implicit_timings_per_n}
\tikzsetnextfilename{implicit_timings_per_n}

\begin{tikzpicture}
\pgfplotsset{every axis/.append style={
    ymajorgrids,
    bar width=5pt,
    width=8cm,
    xtick=data,
    width=8cm, 
    height=5cm, 
    tick label style = {font=\tiny},
    legend style={font=\tiny},
    legend cell align={left},
    legend columns=1,
}
}
\begin{groupplot}[group style={group size=2 by 1, x descriptions at=edge bottom, horizontal sep=10pt}, ]

    \nextgroupplot[xlabel=${n_i}$, title=${(a)}$ ,  ylabel=${\overline{\mathcal{C}}[s]}$, symbolic x coords={4, 8, 12, 16, 24, 32}, xticklabels={$4^2$, $8^2$, $12^2$, $16^2$, $24^2$, $32^2$}, ymin=-1, ymax=30,legend to name=legend_bar_plot_ni, legend image code/.code={
        \draw [#1] (0cm,-0.1cm) rectangle (0.2cm,0.25cm); },]

    \addlegendimage{fill=matred};        
    \addlegendimage{fill=matazure};   
    \addlegendimage{pattern=grid};        
    \addlegendimage{pattern=north east lines};        
    \addlegendimage{pattern=dots};        
    \legend{$\overline{\mathcal{C}}_{\mathrm{assem}}$ , $\overline{\mathcal{C}}_{\mathrm{solve}}$ , HDG-IP, HDG-MX, DG-IP}

    \addplot[ybar, stack plots=y, bar shift = -0.7em, postaction={pattern=grid}, fill=matazure] table [y=solve,x=N, col sep=comma] {sequential_ip.csv};
    \addplot[ybar, stack plots=y, bar shift = -0.7em, postaction={pattern=grid}, fill=matred] table [y=assemble,x=N, col sep=comma] {sequential_ip.csv};
    \addplot[ybar, stack plots=y, bar shift = -0.7em, stack dir =minus, draw=none] table [y=solve,x=N, col sep=comma] {sequential_ip.csv}; 
    \addplot[ybar, stack plots=y, bar shift = -0.7em, stack dir =minus, draw=none] table [y=assemble,x=N, col sep=comma] {sequential_ip.csv};

    \addplot[ybar, stack plots=y, fill=matazure, postaction={pattern=north east lines}] table [y=solve,x=N, col sep=comma] {sequential_mixed.csv};
    \addplot[ybar, stack plots=y, fill=matred, postaction={pattern=north east lines}] table [y=assemble,x=N, col sep=comma] {sequential_mixed.csv};
    \addplot[ybar, stack plots=y, stack dir =minus, draw=none] table [y=solve,x=N, col sep=comma] {sequential_mixed.csv}; 
    \addplot[ybar, stack plots=y, stack dir =minus, draw=none] table [y=assemble,x=N, col sep=comma] {sequential_mixed.csv};

    \addplot[ybar, stack plots=y, bar shift = 0.7em,postaction={pattern=dots},  fill=matazure] table [y=solve,x=N, col sep=comma] {sequential_dg.csv};
    \addplot[ybar, stack plots=y, bar shift = 0.7em,postaction={pattern=dots},  fill=matred] table [y=assemble,x=N, col sep=comma] {sequential_dg.csv};

    \nextgroupplot[xlabel=$n_i$, title=${(b)}$, symbolic x coords={32, 48, 64, 96, 128, 196, 256, 512}, xticklabels={$32^2$, $48^2$, $64^2$, $96^2$, $128^2$, $196^2$, $256^2$, $512^2$}, ymin=-1, ymax=30, yticklabels = {}]

    \addplot[ybar, stack plots=y, bar shift = -0.35em, postaction={pattern=grid}, fill=matazure] table [y=solve,x=N, col sep=comma] {parallel_ip.csv};
    \addplot[ybar, stack plots=y, bar shift = -0.35em, postaction={pattern=grid}, fill=matred] table [y=assemble,x=N, col sep=comma] {parallel_ip.csv};
    \addplot[ybar, stack plots=y, bar shift = -0.35em, stack dir =minus, draw=none] table [y=solve,x=N, col sep=comma] {parallel_ip.csv}; 
    \addplot[ybar, stack plots=y, bar shift = -0.35em, stack dir =minus, draw=none] table [y=assemble,x=N, col sep=comma] {parallel_ip.csv};

    \addplot[ybar, stack plots=y, bar shift = 0.35em, postaction={pattern=north east lines}, fill=matazure] table [y=solve,x=N, col sep=comma] {parallel_mixed.csv};
    \addplot[ybar, stack plots=y, bar shift = 0.35em, postaction={pattern=north east lines}, fill=matred] table [y=assemble,x=N, col sep=comma] {parallel_mixed.csv};

    \coordinate (right_bottom) at (rel axis cs:1,0);
    \coordinate (right_top) at (rel axis cs:1,1);

\end{groupplot}

\path (right_bottom)--(right_top) coordinate[midway] (group center);
\node[inner sep=0pt,xshift=29pt] at(group center) {\pgfplotslegendfromname{legend_bar_plot_ni}};
\end{tikzpicture}
    \caption[Computational Cost of Implicit Evaluation: Assembly and Solve 2]{Computational Cost of the Implicit Solver: Mean wall-clock time (per iteration) in seconds for single-threaded $l=1$ $(a)$ and multi-threaded $l=64$ $(b)$ simulations, both for polynomial order $k=3$.}
    \label{fig::timings_over_n}
\end{figure}
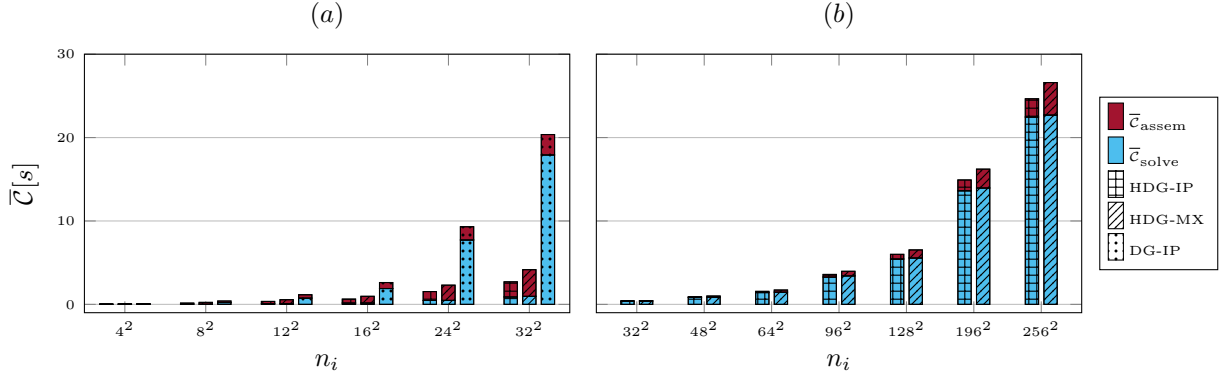

Building on the preceding comparison between HDG-MX and HDG-IP, \Cref{fig::timings_over_n} presents the mean implicit solve and assemble wall-clock time in seconds for polynomial order $k=3$, as a function of the number of implicit elements $n_{i}$, for a single-threaded implementation ($l=1$) and a multi-threaded implementation ($l=64$). Figure (a) underlines the significant performance benefits of HDG formulations over implicit DG schemes, when it comes to solving the linear system. For instance, for $n_{i} = 32^2$ elements, the HDG-IP scheme is approximately $7.5$ times faster, while the HDG-MX scheme is about $5$ times faster than the DG scheme. Figure (b) shows that increasing shared-memory parallelism ($l=64$) yields a pronounced reduction in assembly time—owing to the element-local nature of assembly and local dense operations—while delivering only modest improvements in the solve time. The limited speedup of the solve phase is attributable to the parallel scalability bounds of the chosen sparse direct solver. 

\section{Conclusion}
\label{sec::conclusion}

In this work, we developed an implicit-explicit (IMEX) time integration framework for the compressible Navier-Stokes equations, coupling a DG-IP discretization in the explicit region with two HDG formulations in the implicit region: a mixed formulation (HDG-MX) and a simplified interior-penalty variant (HDG-IP). The framework is designed to address geometry-induced stiffness in highly stretched boundary-layer meshes and is implemented in the open-source \href{https://plederer.github.io/dream_solver/index.html}{\texttt{DreAm}} solver.

Accuracy and implementation correctness were verified using the method of manufactured solutions, confirming the expected temporal convergence rates of the underlying additive Runge-Kutta schemes. Additional validation on an isentropic vortex and laminar NACA-0012 airfoil flow at $\Re_c = 5000$ shows good agreement with reference aerodynamic data.

Computational performance was assessed on a two-dimensional laminar cylinder flow, focusing on the influence of mesh partitioning and geometric stiffness $r_x$. The results show that IMEX schemes can achieve substantial speedups over fully explicit DG, with peak values of $S \approx 50$ for IMEX-IP (DG-IP and HDG-IP) at $r_x \approx 50$. However, performance is strongly dependent on the implicit-explicit partitioning; excessively large implicit regions can offset the gains due to increased nonlinear solve cost.

Among the two HDG variants, HDG-IP consistently outperforms HDG-MX across all tested polynomial orders. This difference is primarily attributed to the higher local assembly cost of HDG-MX, which increases significantly with polynomial degree due to the larger element-level systems, while both formulations exhibit comparable global solve costs after static condensation. At $k=3$ (sequential), HDG-IP reduces the linear-system solve time by up to a factor of $7.5$ compared to a standard implicit DG formulation, and by almost a factor of $5$ compared to HDG-MX, highlighting the efficiency benefit of static condensation and primal formulations.

Overall, these results demonstrate that IMEX schemes can significantly improve efficiency for stiff, boundary-layer-dominated flows, provided that the implicit region is carefully chosen.

Future work will focus on extending the framework to three-dimensional problems, where both stiffness and solver cost increase substantially, and on integrating scalable iterative solvers. Further gains may be achieved by additional static condensation within HDG-MX formulations to reduce assembly overhead, as well as by extending the approach to under-resolved turbulent flows (via WR-LES) where multiple stiffness mechanisms coexist.

\section{Acknowledgments}
\label{sec::acknowledgments}
We would like to thank the Austrian Science Fund (FWF) and the Vienna Scientific Cluster (VSC) for their support.
This research was funded in whole or in part by the Austrian Science Fund (FWF) [10.55776/P35931]. 
Calculations were performed using supercomputer resources provided by the Vienna Scientific Cluster (VSC).

\appendix
\section{Inviscid and Viscous Jacobian Matrices in 2D}
\label{app::jacobians}

\subsection{Inviscid Jacobian Matrices}
\label{sec::inviscid_jacobians}

The inviscid Jacobian matrices in the $x$ and $y$ directions are
\begin{align}
	\vec{A}_{1} 
	=
	\begin{pmatrix}
		0                               &                         1 &                       0 & 1\\
		\overline{\gamma} e_k - u^2     &            (3 - \gamma) u &   - \overline{\gamma} v & \overline{\gamma}\\
		- u v                           &                         v &                       u & 0\\
		u ( \overline{\gamma} e_k - H ) & H - \overline{\gamma} u^2 & - \overline{\gamma} u v & \gamma u
	\end{pmatrix},
	\qquad 
	\vec{A}_{2}
	=
	\begin{pmatrix}
		0                               &                       0 &                         1 & 0\\
		- u v                           &                       v &                         u & 0\\
		\overline{\gamma} e_k - v^2     &   - \overline{\gamma} u &            (3 - \gamma) v & \overline{\gamma}\\
		v ( \overline{\gamma} e_k - H ) & - \overline{\gamma} u v & H - \overline{\gamma} v^2 & \gamma v
	\end{pmatrix},
\end{align}
where $\overline{\gamma} = \gamma - 1$.

\subsection{Viscous Jacobian Matrices}
\label{sec::viscous_jacobians}

The entire \textit{block} diffusion matrix is comprised of
\begin{align}
	\vec{K}
	:=
	\begin{pmatrix}
		\vec{K}_{11} & \vec{K}_{12}\\
		\vec{K}_{21} & \vec{K}_{22}
	\end{pmatrix},
\end{align}
where each of its sub-matrices is defined by
\begin{equation}
	\begin{alignedat}{2}
		\vec{K}_{11}
		&=
		\frac{1}{\rho}
		\begin{pmatrix}
			0                                                      &                               0 &                0 & 0\\
			-\overline{\lambda} u                                  &              \overline{\lambda} &                0 & 0\\
			-\mu v                                                 &                               0 &              \mu & 0\\
			-\mu v^2 - \overline{\lambda} u^2 - \alpha (e_k - e_i) & (\overline{\lambda} - \alpha) u & (\mu - \alpha) v & \alpha
		\end{pmatrix},
		\quad
		\vec{K}_{12}
		&=
		\frac{1}{\rho}
		\begin{pmatrix}
			0                   &   0   &         0 & 0\\
			-\lambda v          &   0   &   \lambda & 0\\
			-\mu u              &   \mu &         0 & 0\\
			-\overline{\mu} u v & \mu v & \lambda u & 0
		\end{pmatrix},\\[2pt]
		\vec{K}_{22}
		&=
		\frac{1}{\rho}
		\begin{pmatrix}
			0                                                 &           0 &                          0 & 0\\
			-\mu u                                            &         \mu &                          0 & 0\\
			-\overline{\lambda} v                             &           0 &         \overline{\lambda} & 0\\
			-\mu u^2 - \overline{\lambda} v^2 - \alpha (e_k - e_i) & (\mu - \alpha) u & (\overline{\lambda} - \alpha) v & \alpha
		\end{pmatrix},
		\quad
		\vec{K}_{21}
		&=
		\frac{1}{\rho}
		\begin{pmatrix}
			0                   &         0 &     0 & 0\\
			-\mu v              &         0 &   \mu & 0\\
			-\lambda u          &   \lambda &     0 & 0\\
			-\overline{\mu} u v & \lambda v & \mu u & 0
		\end{pmatrix},
\end{alignedat}
\end{equation}
where $\alpha = \kappa/c_v$, $\overline{\mu} = \mu + \lambda$ and $\overline{\lambda} = \lambda + 2 \mu$. Recall, $\lambda = -2 \mu / 3$ is the second viscosity.

\section{Physical Boundary Conditions}
\label{app::physical_bcs}

In this work, two types of boundary conditions are considered: farfield boundaries and adiabatic walls. Farfield boundaries are used to truncate the computational domain while allowing outgoing disturbances to leave. At these boundaries, the flow is assumed to approach a prescribed free-stream state, given in conservative variables by
\begin{align} \label{eq::farfield_state}
	\vec{U}_{\infty} 
	=
	(\rho_\infty,\
	\rho_\infty u_\infty,\
	\rho_\infty v_\infty,\
	\rho_\infty E_\infty)^T.
\end{align}

Adiabatic wall boundaries represent impermeable solid surfaces with no heat transfer through the wall. Physically, this corresponds to zero velocity and a vanishing normal heat flux at the wall. These conditions are imposed weakly within the different formulations. The no-slip condition is enforced through a wall state $\vec{U}_w$, defined as
\begin{align} \label{eq::noslip_state}
	\vec{U}_w = (\rho,\, 0,\, 0,\, \rho e_i)^T,
\end{align}
where the velocity components vanish at the wall.

\subsection{Standard Discontinuous Galerkin: Interior Penalty}
\label{sec::sdg_bcs}

\subsubsection{Boundary Condition: Farfield}
\label{subsubsec::sdg_bc_farfield}

The contribution of the farfield boundary, $\Gamma_\infty$, to the weak formulation in \Cref{eq::weak_form_sdg_ip} is defined as
\begin{multline} \label{eq::bilinear_form_sdg_exterior_farfield}
	\mathcal{B}_{\Gamma} 
	:= 
	\mathcal{B}_{\Gamma_{\infty}} \big( \vec{U}_{h}, \vec{U}_{\infty}; \vec{V}_{h} \big)
	=
	\sum_{F \in \Gamma_{\infty}}
	\int\limits_{F}
	\vec{F}^{c}_{n} (\vec{U}_{\Gamma}) \cdot \vec{V}_{h} \, d\bm{s}
	-
	\int\limits_{F}
	\vec{F}^{\nu}_{n} (\vec{U}_{\infty}, \nabla \vec{U}_{h}) \cdot \vec{V}_{h} \, d\bm{s}\\
	-
	\int\limits_{F}
	\vec{K}^T (\vec{U}_{\infty}) \nabla \vec{V}_{h} : \Big( ( \vec{U}_{h} - \vec{U}_{\infty} ) \otimes \vec{n} \Big)\, d\bm{s}
	+
	\int\limits_{F}
	\uptau \vec{K} (\vec{U}_{\infty}) \Big( (\vec{U}_{h} - \vec{U}_{\infty}) \otimes \vec{n} \Big): (\vec{V}_{h} \otimes \vec{n}) \, d\bm{s}.
\end{multline}

The inviscid flux, $\vec{F}^{c}_{n}(\vec{U}_\Gamma)$, is evaluated using a boundary state $\vec{U}_{\Gamma}$, obtained from a characteristic decomposition of the inviscid flux Jacobian in the normal direction, satisfying
\begin{align}
	\mat{A}_{n} (\vec{U}_{h})\, \vec{U}_{\Gamma}
	=
	\mat{A}_{n}^{+} (\vec{U}_{h})\, \vec{U}_{h}
	+
	\mat{A}_{n}^{-} (\vec{U}_{h})\, \vec{U}_{\infty},
\end{align}
where this construction enforces incoming characteristic information from the free-stream state $\vec{U}_{\infty}$, while allowing outgoing characteristics to be determined by the interior solution $\vec{U}_{h}$.

\subsubsection{Boundary Condition: Adiabatic Walls}
\label{subsubsec::sdg_bc_adiabaticwall}

The boundary contribution for adiabatic walls, $\Gamma_{\!w}$, in \Cref{eq::weak_form_sdg_ip}, based on the adjoint-consistent formulation in \cite{hartmann2015generalized}, is given by 
\begin{multline} \label{eq::bilinear_form_sdg_exterior_adiabaticwall}
	\mathcal{B}_{\Gamma}
	:= 
	\mathcal{B}_{\Gamma_{\!w}} \big( \vec{U}_{h}, \vec{U}_{w}; \vec{V}_{h} \big)
	=
	\sum_{F \in \Gamma_{\!w}}
	\int\limits_{F}
	\vec{F}^{c,*}_{n} (\vec{U}_{h}, \vec{U}_{w}) \cdot \vec{V}_{h} \, d\bm{s}
	-
	\int\limits_{F}
	\avg{ \mat{K}_{w} (\vec{U}_{h})} \nabla \vec{U}_{h} : (\vec{V}_{h} \otimes \vec{n}) \, d\bm{s}\\
	-
	\int\limits_{F}
	\mat{K}_{w}^{T} (\vec{U}_{h}) \nabla \vec{V}_{h} : \Big( ( \vec{U}_{h} - \vec{U}_{w} ) \otimes \vec{n} \Big) \, d\bm{s}
	+
	\int\limits_{F}
	\avg{ \mat{K}_{w} (\vec{U}_{h})} \Big( (\vec{U}_{h} - \vec{U}_{w}) \otimes \vec{n} \Big) : (\vec{V}_{h} \otimes \vec{n}) \, d\bm{s},
\end{multline}
where $\vec{U}_{w}$ is the wall state defined in \Cref{eq::noslip_state}.

The inviscid numerical flux $\vec{F}^{c,*}_{n}$ is chosen to satisfy the impermeability constraint at the wall while maintaining consistency and stability. The viscous contribution is constructed using modified viscous matrices $\mat{K}_{w}$, corresponding to an adiabatic wall flux $\vec{F}^{\nu,w}$. These matrices are obtained by enforcing zero heat flux in the normal direction at the wall, i.e. $\vec{q} \cdot \vec{n} = 0$, see \Cref{eq::inviscid_viscous_fluxes,eq::diffusion_block_matrix_product}, such that
\begin{align} \label{eq::adiabatic_viscous_flux_compact}
	\vec{F}^{\nu,w} ( \vec{U}_{h}, \nabla \vec{U}_{h} ) = \mat{K}_{w} (\vec{U}_{h}) \nabla \vec{U}_{h}.
\end{align}

\subsection{Hybridizable Discontinuous Galerkin: Interior Penalty}
\label{sec::hdg_ip_bcs}

\subsubsection{Boundary Condition: Farfield}
\label{subsubsec::hdg_ip_bc_farfield}

The boundary contribution of the farfield boundary, $\Gamma_\infty$, to the weak formulation in \Cref{eq::weak_form_hdg_ip} is enforced through $\widehat{\vec{U}}_{h}$ via the transmissibility conditions using an inviscid flux-vector splitting approach
\begin{align} \label{eq::bilinear_form_hdg_exterior_farfield_facet}
	\widehat{\mathcal{B}}_{\Gamma} 
	:= 
	\widehat{\mathcal{B}}_{\Gamma_{\infty}} \big( \vec{U}_{h}, \widehat{\vec{U}}_{h}, \vec{U}_{\infty}; \vec{V}_{h} \big) 
	=
	\sum_{F \in \Gamma_{\infty}}
	\int\limits_{F}
	\left(
	\mat{A}_{n}^{+} (\widehat{\vec{U}}_{h})\,
	( \widehat{\vec{U}}_{h} - \vec{U}_{h} )
	-
	\mat{A}_{n}^{-} (\widehat{\vec{U}}_{h})\,
	( \widehat{\vec{U}}_{h} - \vec{U}_{\infty} ) 
	\right) \cdot \widehat{\vec{V}}_{h} \, d\bm{s},
\end{align}
 where $\vec{U}_{\infty}$ is defined in \Cref{eq::farfield_state}.

\subsubsection{Boundary Condition: Adiabatic Walls}
\label{subsubsec::hdg_ip_bc_adiabaticwall}

The boundary contribution for adiabatic walls, $\Gamma_{\!w}$, to the weak formulation in \Cref{eq::weak_form_hdg_ip_element} enforces no heat exchange in the normal direction
\begin{multline}
	\label{eq::bilinear_form_hdg_exterior_adiabaticwall_element}
	\mathcal{B}_{\Gamma_{\!w}} 
	:= 
	\mathcal{B}_{\Gamma_{\!w}} \big( \vec{U}_{h}, \widehat{\vec{U}}_{h}; \vec{V}_{h} \big)
	=
	\sum_{F \in \Gamma_{\!w}}
	-
	\int\limits_{F}
	\vec{K}_{w} (\widehat{\vec{U}}_{h}) \nabla \vec{U}_{h} : (\vec{V}_{h} \otimes \vec{n}) \, d\bm{s}\\
	-
	\int\limits_{F}
	\vec{K}_{w}^{T} (\widehat{\vec{U}}_{h}) \nabla \vec{V}_{h} : \widehat{\jump{\vec{U}_{h}}} \, d\bm{s}
	+
	\int\limits_{F}
	\uptau \vec{K}_{w} (\widehat{\vec{U}}_{h}) \widehat{\jump{\vec{U}_{h}}} : (\vec{V}_{h} \otimes \vec{n})\, d\bm{s},
\end{multline}
where the block matrices $\mat{K}_{w}$ fulfill the same relation in \Cref{eq::adiabatic_viscous_flux_compact}.

Additionally, the boundary contribution of $\Gamma_{\!w}$ in \Cref{eq::weak_form_hdg_ip_facet} reduce to a Dirichlet condition enforcing the no-slip constraint. Hence, the transmissibility condition becomes
\begin{align}
	\widehat{\mathcal{B}}_{\Gamma} 
	:= 
	\widehat{\mathcal{B}}_{\Gamma_{\!w}} \big( \widehat{\vec{U}}_{h}, \vec{U}_{w}; \widehat{\vec{V}}_{h} \big)
	=
	\sum_{F \in \Gamma_{\!w}}
	\int\limits_{F}
	\left( \widehat{\vec{U}}_{h} - \vec{U}_w \right) \cdot \widehat{\vec{V}}_{h} \, d\bm{s},
\end{align}
where $\vec{U}_{w}$ is the wall state defined in \Cref{eq::noslip_state}.

\subsection{Hybridizable Discontinuous Galerkin: Mixed Variables}
\label{sec::hdg_mixed_bcs}

\subsubsection{Boundary Condition: Farfield}
\label{subsubsec::hdg_mixed_bc_farfield}

The boundary contribution of the farfield boundary, $\Gamma_\infty$, to the weak formulation in \Cref{eq::weak_form_hdg_mixed_facet}  is enforced in the same fashion as in the HDG-IP formulation, namely by the boundary operator $\widehat{\mathcal{B}}_{\Gamma}$, defined in \Cref{eq::bilinear_form_hdg_exterior_farfield_facet}.

\subsubsection{Boundary Condition: Adiabatic Walls}
\label{subsubsec::hdg_mixed_bc_adiabaticwall}

The adiabatic wall boundary condition in a mixed formulation is enforced through the transmissibility condition, by defining 
\begin{align}
	\widehat{\mathcal{B}}_{\Gamma} = \widehat{\mathcal{B}}_{\Gamma_{\!w}}(\vec{U}_h, \widehat{\vec{U}}_h, \vec{\phi}_h; \widehat{\vec{V}}_h) &:= 
	\sum_{F \in \Gamma_{\!w}} \int_F \left[ \begin{pmatrix}
		\hat{\rho}_h - \rho_h \\[0.8ex] \hat{\rho \VEL}_h \\[0.8ex] \kappa \vec{\phi}_h \cdot \vec{n}
	\end{pmatrix} - \mat{S}_d \begin{pmatrix}
		0 \\ \bm{0} \\ \rho E_h - \hat{\rho E}_h
	\end{pmatrix} \right] \cdot \widehat{\vec{V}}_h \, d\bm{s}.
\end{align}

Unlike HDG-IP (see \Cref{eq::bilinear_form_hdg_exterior_adiabaticwall_element,eq::adiabatic_viscous_flux_compact}), the HDG-MX formulation does not modify the viscous Jacobian matrices to impose adiabatic wall boundary conditions.
\section{IMEX Algorithm}
\label{app::imex_algorithm}

This appendix summarizes the stage-by-stage update used for the IMEX time integration considered in this work. The explicit and implicit components are advanced in a synchronized manner within each Runge-Kutta stage, so that the coupling between the DG and HDG variables is made explicit. For completeness, \Cref{alg::imex_scheme_update} provides the full procedure for one time step, including the treatment of the final update in the presence of the first-same-as-last (FSAL) property for the explicit scheme and stiff accuracy for the implicit scheme.

\begin{algorithm}[h!] 
	\caption{Update procedure for a single time-step using an $s$-stage IMEX scheme.}
	\label{alg::imex_scheme_update}
	\setstretch{1.2}
	\LinesNumbered
	
	\KwIn{DG and HDG solution at time $t=t_{n}$, namely: 
		$\vec{U}_{n}^{ex}$ and $\vec{U}_{n}^{im}$}
	\KwOut{DG and HDG solution at time $t=t_{n+1}$, namely: 
		$\vec{U}_{n+1}^{ex}$ and $\vec{U}_{n+1}^{im}$}
	\KwReq{synchronized RK-IMEX schemes: $\overline{c}_{i+1} = c_{i}$ for $i = 1, \cdots, s$.}
	
	\textbf{Initial step:} set initial stage values: 
	$\vec{y}_{0}^{im} = \vec{U}_{n}^{im}$ and $\vec{y}_{1}^{ex} = \vec{U}_{n}^{ex}$\\
	
	\tcp{loop over each stage and update the stage-solutions.}
	\For{$i=1$ \KwTo $s$}{
		\tcp{update the explicit solution first, then the implicit one.}
		\textbf{Explicit stage:} compute $\vec{y}_{i+1}^{ex}$ using \Cref{eq::erk_stage_update}
		
		\textbf{Implicit stage:} compute $\vec{y}_{i}^{im}$ by solving the nonlinear system in \Cref{eq::dirk_stage_update_solution,eq::dirk_stage_update_facet,eq::dirk_stage_update_mixed}
	}
	
	\tcp{the final explicit solution depends on first same as last (FSAL) property.}
	\uIf{explicit scheme does not have the FSAL property (i.e. $\overline{b}_{i} \neq \overline{a}_{s+1,i}$ or $c_{s+1} \neq 1$)}{
		\textbf{Final explicit step:} compute $\vec{U}_{n+1}^{ex}$ using \Cref{eq::erk_solution_update}
	}
	\Else{
		\tcp{the last stage corresponds to the solution, since $\overline{b}_{i} = \overline{a}_{s+1,i}$ and $c_{s+1} = 1$.}
		\textbf{FSAL case:} set $\vec{U}_{n+1}^{ex} = \vec{y}_{s+1}^{ex}$
	}
	
	\tcp{the final implicit solution depends on whether the scheme is stiffly accurate.}
	\uIf{implicit scheme is not stiffly accurate (i.e. $b_{i} \neq a_{si}$ or $c_{s} \neq 1$)}{
		\textbf{Final implicit step:} compute $\vec{U}_{n+1}^{im}$ using \Cref{eq::dirk_solution_update}
	}
	\Else{
		\tcp{the last stage corresponds to the solution, since $b_{i} = a_{si}$ and $c_{s} = 1$.}
		\textbf{Stiffly accurate case:} set $\vec{U}_{n+1}^{im} = \vec{y}_{s}^{im}$
	}
\end{algorithm}

\section{Complementary Numerical Results}
\label{app::complementary_results}
\pgfplotsset{table/search path={./figures/cylinder_data/}}

This appendix collects complementary numerical results that support the observations made in \Cref{sec::experiments}-

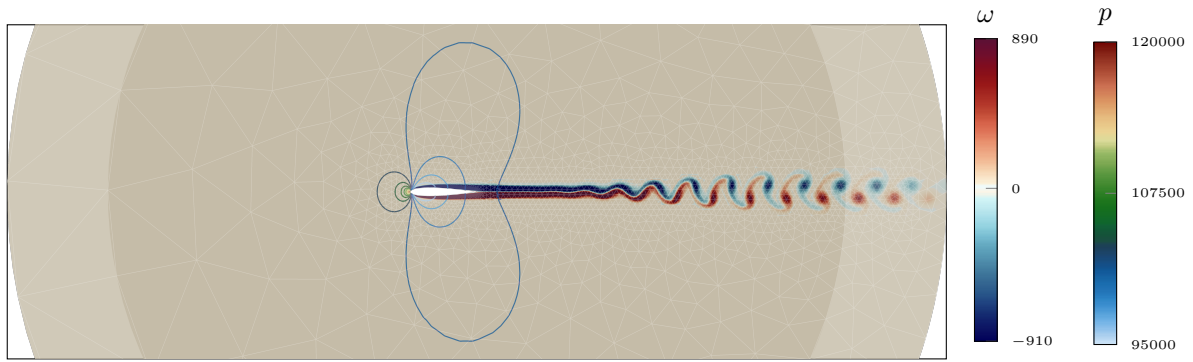
\begin{figure}[htb]
    \centering

    \tikzsetnextfilename{naca_contours}
    \begin{tikzpicture}[>=latex]
        \pgfplotsset{table/search path={./figures/naca_data/}}
        \graphicspath{{./figures/naca_data/}}
        \pgfplotsset{every axis/.append style={
                    width=14cm,
                    height=6cm,
                    thick,
                    axis line style={thin},
                    tick label style = {font=\tiny},
                    ytick pos=left,},
            xtick={\empty},
            ytick={\empty},
            legend style={font=\tiny},
            legend cell align={left},
            legend columns=1,
            transpose legend,
colormap={pressure_color}{
    rgb255(0.0cm)=(212,232,249);
    rgb255(0.0125cm)=(191,219,244);
    rgb255(0.025cm)=(177,211,239);
    rgb255(0.05cm)=(144,193,235);
    rgb255(0.075cm)=(114,179,230);
    rgb255(0.1cm)=(87,163,219);
    rgb255(0.125cm)=(63,145,209);
    rgb255(0.15cm)=(45,133,198);
    rgb255(0.16cm)=(38,124,191);
    rgb255(0.18cm)=(33,114,181);
    rgb255(0.2cm)=(26,109,175);
    rgb255(0.21cm)=(24,103,168);
    rgb255(0.22cm)=(22,100,163);
    rgb255(0.23cm)=(21,93,158);
    rgb255(0.24cm)=(17,89,153);
    rgb255(0.25cm)=(17,84,144);
    rgb255(0.26cm)=(19,79,137);
    rgb255(0.27cm)=(22,77,130);
    rgb255(0.28cm)=(24,73,121);
    rgb255(0.29cm)=(26,70,114);
    rgb255(0.3cm)=(28,68,105);
    rgb255(0.31cm)=(28,66,96);
    rgb255(0.32cm)=(28,63,89);
    rgb255(0.33cm)=(28,68,82);
    rgb255(0.34cm)=(26,77,66);
    rgb255(0.35cm)=(24,79,61);
    rgb255(0.36cm)=(21,82,58);
    rgb255(0.37cm)=(19,86,56);
    rgb255(0.38cm)=(17,91,54);
    rgb255(0.39cm)=(17,96,54);
    rgb255(0.4cm)=(15,102,45);
    rgb255(0.425cm)=(19,107,36);
    rgb255(0.45cm)=(22,112,29);
    rgb255(0.475cm)=(31,119,29);
    rgb255(0.5cm)=(47,128,38);
    rgb255(0.525cm)=(65,137,47);
    rgb255(0.55cm)=(82,147,59);
    rgb255(0.575cm)=(103,158,70);
    rgb255(0.6cm)=(128,170,84);
    rgb255(0.63cm)=(151,186,102);
    rgb255(0.65cm)=(188,200,124);
    rgb255(0.67cm)=(219,219,154);
    rgb255(0.7cm)=(235,212,147);
    rgb255(0.75cm)=(230,186,126);
    rgb255(0.8cm)=(219,149,98);
    rgb255(0.85cm)=(204,112,82);
    rgb255(0.9cm)=(172,75,52);
    rgb255(0.95cm)=(140,42,28);
    rgb255(0.975cm)=(121,21,12);
    rgb255(1.0cm)=(114,1,0);
},
colormap={vorticity_color}{
    rgb255(0.0cm)=(0,0,89);
    rgb255(0.03125cm)=(10,15,96);
    rgb255(0.0625cm)=(15,29,105);
    rgb255(0.09375cm)=(22,47,114);
    rgb255(0.125cm)=(31,66,128);
    rgb255(0.15625cm)=(40,86,137);
    rgb255(0.1875cm)=(51,100,144);
    rgb255(0.21875cm)=(61,116,153);
    rgb255(0.25cm)=(73,133,165);
    rgb255(0.28125cm)=(86,151,179);
    rgb255(0.3125cm)=(98,167,191);
    rgb255(0.34375cm)=(119,188,209);
    rgb255(0.375cm)=(145,209,223);
    rgb255(0.40625cm)=(167,221,232);
    rgb255(0.4375cm)=(191,233,239);
    rgb255(0.46875cm)=(209,244,246);
    rgb255(0.49cm)=(251,244,230);
    rgb255(0.51cm)=(239,251,251);
    rgb255(0.52cm)=(251,240,216);
    rgb255(0.54cm)=(249,228,200);
    rgb255(0.5625cm)=(246,212,177);
    rgb255(0.59375cm)=(242,186,149);
    rgb255(0.625cm)=(237,165,130);
    rgb255(0.65625cm)=(232,144,110);
    rgb255(0.6875cm)=(223,117,89);
    rgb255(0.71875cm)=(214,98,73);
    rgb255(0.75cm)=(193,75,54);
    rgb255(0.78125cm)=(179,54,42);
    rgb255(0.8125cm)=(165,40,33);
    rgb255(0.84375cm)=(153,24,24);
    rgb255(0.875cm)=(140,17,24);
    rgb255(0.90625cm)=(128,12,31);
    rgb255(0.9375cm)=(114,14,43);
    rgb255(0.96875cm)=(102,14,49);
    rgb255(1.0cm)=(89,17,54);},
        }

        \begin{groupplot}[group style={group size=1 by 2, vertical sep=0.5cm, horizontal sep=0.5cm}]

            \nextgroupplot[xmin=-6, xmax=8, ymin=-2, ymax=2]

            \addplot graphics[xmin=-6,ymin=-3,xmax=8,ymax=3] {instantaneous.png};


            \coordinate (top) at (rel axis cs:1,1);
            \coordinate (bottom) at (rel axis cs:1,0);

        \end{groupplot}

        \begin{axis}[%
                hide axis,
                scale only axis,
                colorbar=true,
                xshift=0cm,
                colorbar to name=vorticity_colorbar,
                colormap name=vorticity_color,
                colorbar style={title=$\omega$,
                        height=4cm,
                        width=0.3cm,
                        ytick = {-910, 0, 890},
                        yticklabels = {$-910$, $0$, $890$},
                        ytick pos=right,
                        scaled y ticks=false,
                        yticklabel style={anchor=west, xshift=0.5em},
                        title style={align=center, yshift=0.0em},
                    },
                point meta min = -910,
                point meta max = 890
            ]
            \addplot [draw=none] coordinates {(0,0)};
        \end{axis}

        \begin{axis}[%
                hide axis,
                scale only axis,
                colorbar=true,
                xshift=0cm,
                colorbar to name=pressure_colorbar,
                colormap name=pressure_color,
                colorbar style={title=$p$,
                        height=4cm,
                        width=0.3cm,
                        ytick = {95000, 107500, 120000},
                        yticklabels = {$95000$, $107500$, $120000$},
                        ytick pos=right,
                        scaled y ticks=false,
                        yticklabel style={anchor=west, xshift=0.5em},
                        title style={align=center, yshift=0.0em},
                    },
                point meta min = 95000,
                point meta max = 120000
            ]
            \addplot [draw=none] coordinates {(0,0)};
        \end{axis}

        \path (top)--(bottom) coordinate[midway] (group center);
        \node[right=1em,inner sep=0pt, yshift=0.5em] at(group center -| current bounding box.east) {\pgfplotscolorbarfromname{vorticity_colorbar}};

        \path (top)--(bottom) coordinate[midway] (group center);
        \node[right=1.5em,inner sep=0pt, yshift=0.5em] at(group center -| current bounding box.east) {\pgfplotscolorbarfromname{pressure_colorbar}};


    \end{tikzpicture}

    \caption[Laminar Flow over a NACA-0012 Airfoil: Vorticity and Pressure]{Laminar Flow over a NACA-0012 Airfoil: Instantaneous surface plot of the vorticity $\omega\, [\textrm{m}/\textrm{s}^2]$ and contour plot of the pressure $p\, [\textrm{Pa}]$.}
    \label{fig::naca0012_contours}
\end{figure}

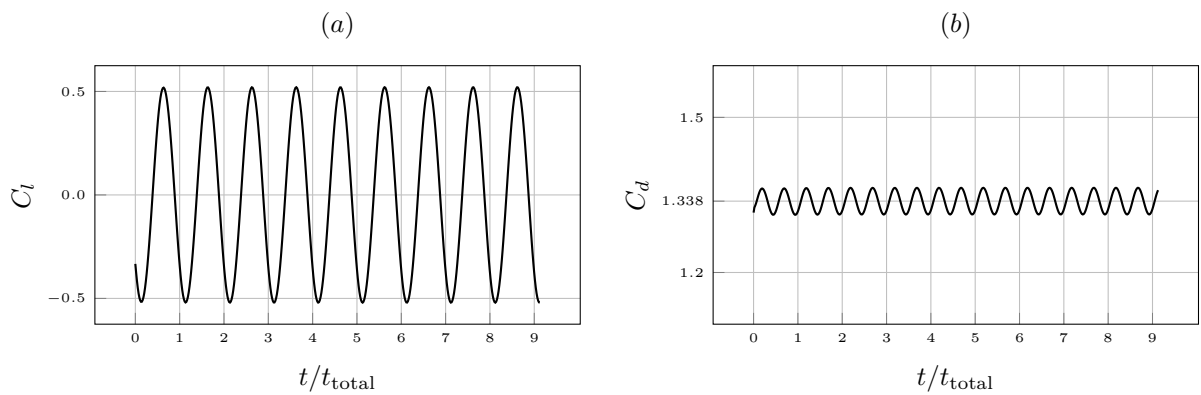
\begin{figure}[htb]
    \centering
    \tikzsetnextfilename{cylinder_coefficients}
    \begin{tikzpicture}[>=latex]

        \pgfplotsset{every axis/.append style={
                    width=8cm,
                    height=5cm,
                    mark size = 1pt,
                    mark repeat = 20,
                    xtick pos=bottom,
                    xmajorgrids=true, ymajorgrids=true,
                    axis line style={thin},
                    tick label style = {font=\tiny},
                    ytick pos=left,},
        }

        \begin{groupplot}[group style={group size=2 by 1, horizontal sep=50pt, vertical sep=0.5cm}]

            \nextgroupplot[ylabel={$C_l$}, xlabel=$t/t_{\text{total}}$, xtick={0.        , 0.10958904, 0.21917808, 0.32876712, 0.43835616,
                    0.54794521, 0.65753425, 0.76712329, 0.87671233, 0.98630137},xticklabels={$0$, $1$, $2$, $3$, $4$, $5$, $6$, $7$, $8$, $9$}, ytick = {-0.5, 0.0, 0.5}, yticklabels = {$-0.5$, $0.0$, $0.5$}, ylabel shift = -2pt, xlabel shift = 2pt]
            \addplot[thick]  table [y=cl, x=t, col sep=comma]{coefficients.dat};
            \path (rel axis cs:0,1)--(rel axis cs:1,1) coordinate[midway] (group center_left);

            \nextgroupplot[ylabel={$C_d$}, xlabel=$t/t_{\text{total}}$, xtick={0.        , 0.10958904, 0.21917808, 0.32876712, 0.43835616,
                    0.54794521, 0.65753425, 0.76712329, 0.87671233, 0.98630137},xticklabels={$0$, $1$, $2$, $3$, $4$, $5$, $6$, $7$, $8$, $9$}, ytick = {1.2, 1.338, 1.5}, yticklabels = {$1.2$, $1.338$, $1.5$}, ymin=1.1, ymax=1.6, ylabel shift = -2pt, xlabel shift = 2pt]
            \addplot[thick]  table [y=cd, x=t, col sep=comma]{coefficients.dat};
            \path (rel axis cs:0,1)--(rel axis cs:1,1) coordinate[midway] (group center_right);

        \end{groupplot}

        \node[inner sep=0pt, yshift=15pt] at(group center_left) {$(a)$};
        \node[inner sep=0pt, yshift=15pt] at(group center_right) {$(b)$};




    \end{tikzpicture}

    \caption[Speedup - Laminar flow over a Cylinder: Lift and Drag]{Speedup - Laminar flow over a Cylinder: Lift coefficient ($a$) and drag coefficient ($b$) over multiple shedding periods~$t_{\text{total}}$.}
    \label{fig::cylinder_coeffcients}
\end{figure}

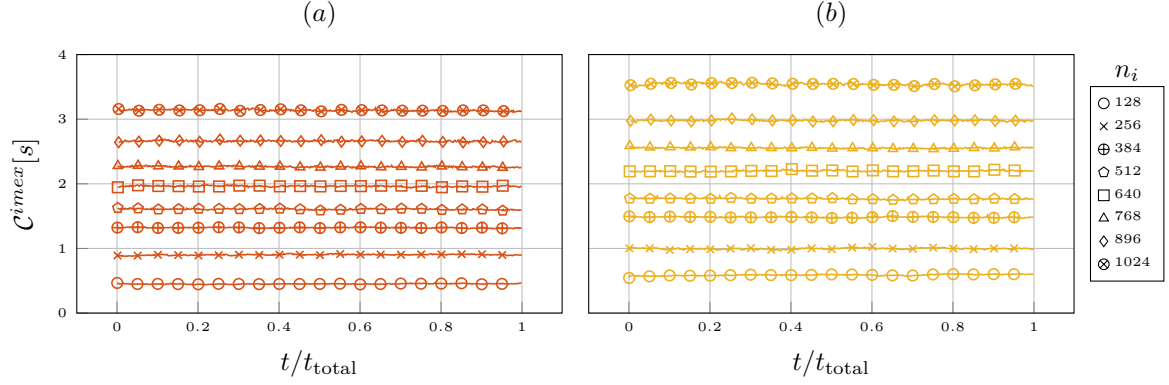
\begin{figure}[htb]
    \centering
    \tikzsetnextfilename{cylinder_domain_runtimes}
    \begin{tikzpicture}[>=latex]

        \pgfplotsset{every axis/.append style={
                    width=8cm,
                    height=5cm,
                    mark size = 2pt,
                    mark repeat = 25,
                    thick,
                    xtick pos=bottom,
                    xmajorgrids=true, ymajorgrids=true,
                    axis line style={thin},
                    tick label style = {font=\tiny},
                    ytick pos=left,},
            legend style={font=\tiny},
            legend cell align={left},
            legend columns=1,
            transpose legend,
        }

        \begin{groupplot}[group style={group size=2 by 1, horizontal sep=10pt, vertical sep=0.5cm}]

            \nextgroupplot[ylabel={$\mathcal{C}^{imex}[s]$}, xlabel={$t/t_{\text{total}}$}, ymin=0, ymax=4]

            \addplot[semithick, matorange, mark=o]  table [y=mean_imp, x=t, col sep=comma]{mean_runtimes_ip_4.csv};
            \addplot[semithick, matorange, mark=x]  table [y=mean_imp, x=t, col sep=comma]{mean_runtimes_ip_8.csv};
            \addplot[semithick, matorange, mark=oplus]  table [y=mean_imp, x=t, col sep=comma]{mean_runtimes_ip_12.csv};
            \addplot[semithick, matorange, mark=pentagon]  table [y=mean_imp, x=t, col sep=comma]{mean_runtimes_ip_16.csv};
            \addplot[semithick, matorange, mark=square]  table [y=mean_imp, x=t, col sep=comma]{mean_runtimes_ip_20.csv};
            \addplot[semithick, matorange, mark=triangle]  table [y=mean_imp, x=t, col sep=comma]{mean_runtimes_ip_24.csv};
            \addplot[semithick, matorange, mark=diamond]  table [y=mean_imp, x=t, col sep=comma]{mean_runtimes_ip_28.csv};
            \addplot[semithick, matorange, mark=otimes]  table [y=mean_imp, x=t, col sep=comma]{mean_runtimes_ip_32.csv};

            \path (rel axis cs:0,1)--(rel axis cs:1,1) coordinate[midway] (group center_left);

            \nextgroupplot[xlabel={$t/t_{\text{total}}$}, y tick style={draw=none}, yticklabels={}, ymin=0, ymax=4]

            \addplot[semithick, matyellow, mark=o]  table [y=mean_imp, x=t, col sep=comma]{mean_runtimes_mixed_4.csv};
            \addplot[semithick, matyellow, mark=x]  table [y=mean_imp, x=t, col sep=comma]{mean_runtimes_mixed_8.csv};
            \addplot[semithick, matyellow, mark=oplus]  table [y=mean_imp, x=t, col sep=comma]{mean_runtimes_mixed_12.csv};
            \addplot[semithick, matyellow, mark=pentagon]  table [y=mean_imp, x=t, col sep=comma]{mean_runtimes_mixed_16.csv};
            \addplot[semithick, matyellow, mark=square]  table [y=mean_imp, x=t, col sep=comma]{mean_runtimes_mixed_20.csv};
            \addplot[semithick, matyellow, mark=triangle]  table [y=mean_imp, x=t, col sep=comma]{mean_runtimes_mixed_24.csv};
            \addplot[semithick, matyellow, mark=diamond]  table [y=mean_imp, x=t, col sep=comma]{mean_runtimes_mixed_28.csv};
            \addplot[semithick, matyellow, mark=otimes]  table [y=mean_imp, x=t, col sep=comma]{mean_runtimes_mixed_32.csv};

            \coordinate (right_bottom) at (rel axis cs:1,0);
            \coordinate (right_top) at (rel axis cs:1,1);
            \path (rel axis cs:0,1)--(rel axis cs:1,1) coordinate[midway] (group center_right);

        \end{groupplot}
        \path (right_bottom)--(right_top) coordinate[midway] (group center);
        \node[inner sep=0pt,xshift=20pt] at(group center) {\pgfplotslegendfromname{legend_speedup_linear}};
        \node[inner sep=0pt,xshift=20pt, yshift=-7pt] at(right_top) {$n_i$};
        \node[inner sep=0pt, yshift=15pt] at(group center_left) {$(a)$};
        \node[inner sep=0pt, yshift=15pt] at(group center_right) {$(b)$};

    \end{tikzpicture}

    \caption[Speedup - Laminar flow over a Cylinder: Runtimes]{Speedup - Laminar flow over a Cylinder: IMEX wall clock time per time step for IMEX-IP $(a)$ and IMEX-MX $(b)$
        over a shedding period $t_{\text{total}}$ for varying number of implicit elements $n_i$.}
    \label{fig::cylinder_domain_runtimes}
\end{figure}

\begin{figure}[htb]
    \centering
    \tikzsetnextfilename{pulse_domain_runtimes}
    \begin{tikzpicture}[>=latex]

        \pgfplotsset{every axis/.append style={
                    width=8cm,
                    height=5cm,
                    mark size = 2pt,
                    mark repeat = 25,
                    thick,
                    xtick pos=bottom,
                    xmajorgrids=true, ymajorgrids=true,
                    axis line style={thin},
                    tick label style = {font=\tiny},
                    ytick pos=left,},
            legend style={font=\tiny},
            legend cell align={left},
            legend columns=1,
            transpose legend,
            table/search path={./figures/pulse_data/}
        }

        \begin{groupplot}[group style={group size=2 by 1, horizontal sep=40pt, vertical sep=0.5cm}]

            \nextgroupplot[ylabel={$\overline{\mathcal{C}}_{\mathrm{assem}}^{\mathrm{IP}}[s]$}, xlabel={$m$}, ymin=-0.1, ymax=3.1, ylabel shift = -2pt]

            \addplot[semithick, matred]  table [y=assemble, x=m, col sep=comma]   {1x64_ip_solver.csv};
            \addplot[semithick, matazure]  table [y=assemble, x=m, col sep=comma] {2x64_ip_solver.csv};
            \addplot[semithick, matblue]  table [y=assemble, x=m, col sep=comma]  {3x64_ip_solver.csv};
            \addplot[semithick, matgreen]  table [y=assemble, x=m, col sep=comma] {4x64_ip_solver.csv};
            \addplot[semithick, matpurple]  table [y=assemble, x=m, col sep=comma]{5x64_ip_solver.csv};

            \path (rel axis cs:0,1)--(rel axis cs:1,1) coordinate[midway] (group center_left);

            \nextgroupplot[ylabel={$\overline{\mathcal{C}}_{\mathrm{solve}}^{\mathrm{IP}}[s]$}, xlabel={$m$}, , ymin=-0.5, ymax=15.5, ylabel shift = -2pt]

            \addplot[semithick, matred]  table [y=solve, x=m, col sep=comma]   {1x64_ip_solver.csv};
            \addplot[semithick, matazure]  table [y=solve, x=m, col sep=comma] {2x64_ip_solver.csv};
            \addplot[semithick, matblue]  table [y=solve, x=m, col sep=comma]  {3x64_ip_solver.csv};
            \addplot[semithick, matgreen]  table [y=solve, x=m, col sep=comma] {4x64_ip_solver.csv};
            \addplot[semithick, matpurple]  table [y=solve, x=m, col sep=comma]{5x64_ip_solver.csv};

            \coordinate (right_bottom) at (rel axis cs:1,0);
            \coordinate (right_top) at (rel axis cs:1,1);
            \path (rel axis cs:0,1)--(rel axis cs:1,1) coordinate[midway] (group center_right);

        \end{groupplot}
        \path (right_bottom)--(right_top) coordinate[midway] (group center);
        \node[inner sep=0pt,xshift=15pt] at(group center) {\pgfplotslegendfromname{legend_bar_plot_k}};
        \node[inner sep=0pt,xshift=15pt, yshift=-7pt] at(right_top) {$k$};
        \node[inner sep=0pt, yshift=15pt] at(group center_left) {$(a)$};
        \node[inner sep=0pt, yshift=15pt] at(group center_right) {$(b)$};

    \end{tikzpicture}

    \caption[Computational Cost of the Implicit Solver: Runtimes]{Computational Cost of the Implicit Solver: IMEX-IP assembly $(a)$ and solve $(b)$ wall clock time for $m=500$ iterations, $n_i = 128^2$ implicit elements, in a multi-threaded~$l=64$ environment for varying polynomial order~$k$.}
    \label{fig::pulse_domain_runtimes}
\end{figure}
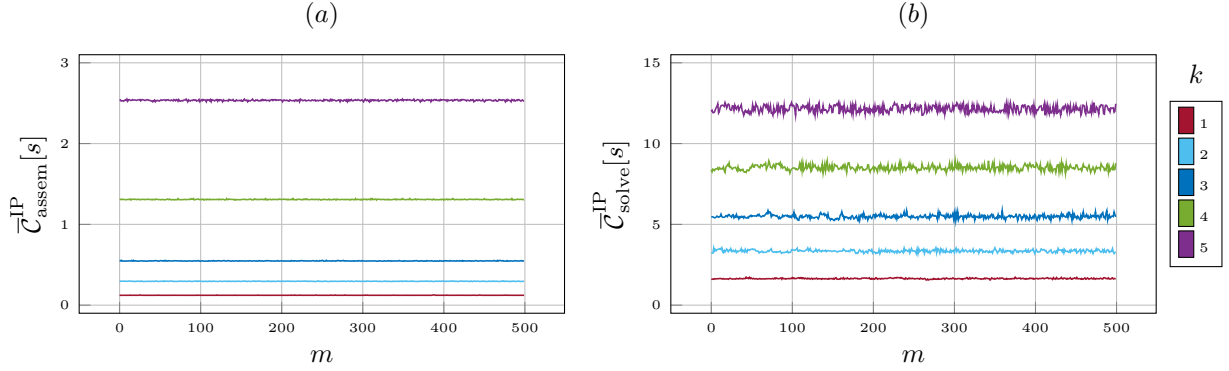

\bibliography{literature}
\end{document}